\documentclass[12pt]{article}
\usepackage{amsmath,latexsym,amssymb}
\begin{document}
\headsep 0.5 true cm

\begin{center}
{\Large\bf Modular curves, invariant theory and $E_8$}
\vskip 1.0 cm
{\large\bf Lei Yang}
\end{center}
\vskip 1.5 cm

\begin{center}
{\large\bf Abstract}
\end{center}

\vskip 0.5 cm

  The $E_8$ root lattice can be constructed from the modular
curve $X(13)$ by the invariant theory for the simple group
$\text{PSL}(2, 13)$. This gives a different construction of the
$E_8$ root lattice. It also gives an explicit construction of
the modular curve $X(13)$.

\begin{center}
{\large\bf Contents}
\end{center}
$$\aligned
  &1. \text{\quad Introduction}\\
  &2. \text{\quad Standard construction: from the icosahedron to $E_8$}\\
  &3. \text{\quad Modular curve $X(13)$ and invariant theory for $\text{PSL}(2, 13)$}\\
  &4. \text{\quad A different construction: from the modular curve $X(13)$ to $E_8$}\\
  &5. \text{\quad An explicit construction of the modular curve $X(13)$}
\endaligned$$

\begin{center}
{\large\bf 1. Introduction}
\end{center}

  The $E_8$ root lattice (see \cite{G}) occurs in: the equation of the
$E_8$-singularity (theory of singularities), the Barlow surface, which
is homeomorphic but not diffeomorphic to $\mathbb{CP}^2 \# 8 \overline{\mathbb{CP}}^2$
(differential topology of $4$-manifolds), and the automorphism group of
the configuration of $120$ tritangent planes of Bring's curve
(representation theory and classical algebraic geometry) (see section
two for more details).

  In the present paper, we will give a different relation which connects
$E_8$ with the modular curve $X(13)$. We will show that the $E_8$ root
lattice can be constructed from the modular curve $X(13)$ by the invariant 
theory for the simple group $\text{PSL}(2, 13)$.

  Let us begin with the invariant theory for $\text{PSL}(2, 13)$. Recall
that the six-dimensional representation (the Weil representation) of the
finite simple group $\text{PSL}(2, 13)$ of order $1092$, which acts on the
five-dimensional projective space
$\mathbb{P}^5=\{ (z_1, z_2, z_3, z_4, z_5, z_6): z_i \in \mathbb{C} \quad
(i=1, 2, 3, 4, 5, 6) \}$. This representation is defined over the cyclotomic
field $\mathbb{Q}(e^{\frac{2 \pi i}{13}})$. Put
$$S=-\frac{1}{\sqrt{13}} \begin{pmatrix}
  \zeta^{12}-\zeta & \zeta^{10}-\zeta^3 & \zeta^4-\zeta^9
& \zeta^5-\zeta^8 & \zeta^2-\zeta^{11} & \zeta^6-\zeta^7\\
  \zeta^{10}-\zeta^3 & \zeta^4-\zeta^9 & \zeta^{12}-\zeta
& \zeta^2-\zeta^{11} & \zeta^6-\zeta^7 & \zeta^5-\zeta^8\\
  \zeta^4-\zeta^9 & \zeta^{12}-\zeta & \zeta^{10}-\zeta^3
& \zeta^6-\zeta^7 & \zeta^5-\zeta^8 & \zeta^2-\zeta^{11}\\
  \zeta^5-\zeta^8 & \zeta^2-\zeta^{11} & \zeta^6-\zeta^7
& \zeta-\zeta^{12} & \zeta^3-\zeta^{10} & \zeta^9-\zeta^4\\
  \zeta^2-\zeta^{11} & \zeta^6-\zeta^7 & \zeta^5-\zeta^8
& \zeta^3-\zeta^{10} & \zeta^9-\zeta^4 & \zeta-\zeta^{12}\\
  \zeta^6-\zeta^7 & \zeta^5-\zeta^8 & \zeta^2-\zeta^{11}
& \zeta^9-\zeta^4 & \zeta-\zeta^{12} & \zeta^3-\zeta^{10}
\end{pmatrix} $$
and
$$T=\text{diag}(\zeta^7, \zeta^{11}, \zeta^8, \zeta^6, \zeta^2, \zeta^5)$$
where $\zeta=\exp(2 \pi i/13)$. We have
$$S^2=T^{13}=(ST)^3=1.\eqno{(1.1)}$$
Let
$G=\langle S, T \rangle$, then $G \cong \text{PSL}(2, 13)$.
We construct some $G$-invariant polynomials in six variables
$z_1, \ldots, z_6$. Let
$$w_{\infty}=13 \mathbf{A}_0^2, \quad
  w_{\nu}=(\mathbf{A}_0+\zeta^{\nu} \mathbf{A}_1+\zeta^{4 \nu} \mathbf{A}_2+\zeta^{9 \nu}
  \mathbf{A}_3+\zeta^{3 \nu} \mathbf{A}_4+\zeta^{12 \nu} \mathbf{A}_5+\zeta^{10 \nu}
  \mathbf{A}_6)^2\eqno{(1.2)}$$
for $\nu=0, 1, \ldots, 12$, where the senary quadratic forms (quadratic forms
in six variables) $\mathbf{A}_j$ $(j=0, 1, \ldots, 6)$ are given by
$$\left\{\aligned
  \mathbf{A}_0 &=z_1 z_4+z_2 z_5+z_3 z_6,\\
  \mathbf{A}_1 &=z_1^2-2 z_3 z_4,\\
  \mathbf{A}_2 &=-z_5^2-2 z_2 z_4,\\
  \mathbf{A}_3 &=z_2^2-2 z_1 z_5,\\
  \mathbf{A}_4 &=z_3^2-2 z_2 z_6,\\
  \mathbf{A}_5 &=-z_4^2-2 z_1 z_6,\\
  \mathbf{A}_6 &=-z_6^2-2 z_3 z_5.
\endaligned\right.\eqno{(1.3)}$$
Then $w_{\infty}$, $w_{\nu}$ for $\nu=0, \ldots, 12$ are the roots of a
polynomial of degree fourteen. The corresponding equation is just the
Jacobian equation of degree fourteen (see \cite{K}, pp.161-162). On the
other hand, set
$$\delta_{\infty}=13^2 \mathbf{G}_0, \quad
  \delta_{\nu}=-13 \mathbf{G}_0+\zeta^{\nu} \mathbf{G}_1+\zeta^{2 \nu}
                \mathbf{G}_2+\cdots+\zeta^{12 \nu} \mathbf{G}_{12}\eqno{(1.4)}$$
for $\nu=0, 1, \ldots, 12$, where the senary sextic forms (i.e., sextic forms
in six variables) $\mathbf{G}_j$ $(j=0, 1, \ldots, 12)$ are given by
$$\left\{\aligned
  \mathbf{G}_0 =&\mathbf{D}_0^2+\mathbf{D}_{\infty}^2,\\
  \mathbf{G}_1 =&-\mathbf{D}_7^2+2 \mathbf{D}_0 \mathbf{D}_1+10 \mathbf{D}_{\infty}
                 \mathbf{D}_1+2 \mathbf{D}_2 \mathbf{D}_{12}+\\
                &-2 \mathbf{D}_3 \mathbf{D}_{11}-4 \mathbf{D}_4 \mathbf{D}_{10}
                 -2 \mathbf{D}_9 \mathbf{D}_5,\\
  \mathbf{G}_2 =&-2 \mathbf{D}_1^2-4 \mathbf{D}_0 \mathbf{D}_2+6 \mathbf{D}_{\infty}
                 \mathbf{D}_2-2 \mathbf{D}_4 \mathbf{D}_{11}+\\
                &+2 \mathbf{D}_5 \mathbf{D}_{10}-2 \mathbf{D}_6 \mathbf{D}_9-2
                 \mathbf{D}_7 \mathbf{D}_8,\\
  \mathbf{G}_3 =&-\mathbf{D}_8^2+2 \mathbf{D}_0 \mathbf{D}_3+10 \mathbf{D}_{\infty}
                 \mathbf{D}_3+2 \mathbf{D}_6 \mathbf{D}_{10}+\\
                &-2 \mathbf{D}_9 \mathbf{D}_7-4 \mathbf{D}_{12} \mathbf{D}_4
                 -2 \mathbf{D}_1 \mathbf{D}_2,\\
  \mathbf{G}_4 =&-\mathbf{D}_2^2+10 \mathbf{D}_0 \mathbf{D}_4-2 \mathbf{D}_{\infty}
                 \mathbf{D}_4+2 \mathbf{D}_5 \mathbf{D}_{12}+\\
                &-2 \mathbf{D}_9 \mathbf{D}_8-4 \mathbf{D}_1 \mathbf{D}_3-2
                 \mathbf{D}_{10} \mathbf{D}_7,\\
  \mathbf{G}_5 =&-2 \mathbf{D}_9^2-4 \mathbf{D}_0 \mathbf{D}_5+6 \mathbf{D}_{\infty}
                 \mathbf{D}_5-2 \mathbf{D}_{10} \mathbf{D}_8+\\
                &+2 \mathbf{D}_6 \mathbf{D}_{12}-2 \mathbf{D}_2 \mathbf{D}_3
                 -2 \mathbf{D}_{11} \mathbf{D}_7,\\
  \mathbf{G}_6 =&-2 \mathbf{D}_3^2-4 \mathbf{D}_0 \mathbf{D}_6+6 \mathbf{D}_{\infty}
                 \mathbf{D}_6-2 \mathbf{D}_{12} \mathbf{D}_7+\\
                &+2 \mathbf{D}_2 \mathbf{D}_4-2 \mathbf{D}_5 \mathbf{D}_1-2
                 \mathbf{D}_8 \mathbf{D}_{11},\\
  \mathbf{G}_7 =&-2 \mathbf{D}_{10}^2+6 \mathbf{D}_0 \mathbf{D}_7+4 \mathbf{D}_{\infty}
                 \mathbf{D}_7-2 \mathbf{D}_1 \mathbf{D}_6+\\
                &-2 \mathbf{D}_2 \mathbf{D}_5-2 \mathbf{D}_8 \mathbf{D}_{12}-2
                 \mathbf{D}_9 \mathbf{D}_{11},\\
  \mathbf{G}_8 =&-2 \mathbf{D}_4^2+6 \mathbf{D}_0 \mathbf{D}_8+4 \mathbf{D}_{\infty}
                 \mathbf{D}_8-2 \mathbf{D}_3 \mathbf{D}_5+\\
                &-2 \mathbf{D}_6 \mathbf{D}_2-2 \mathbf{D}_{11} \mathbf{D}_{10}-2
                 \mathbf{D}_1 \mathbf{D}_7,\\
  \mathbf{G}_9 =&-\mathbf{D}_{11}^2+2 \mathbf{D}_0 \mathbf{D}_9+10 \mathbf{D}_{\infty}
                 \mathbf{D}_9+2 \mathbf{D}_5 \mathbf{D}_4+\\
                &-2 \mathbf{D}_1 \mathbf{D}_8-4 \mathbf{D}_{10} \mathbf{D}_{12}-2
                 \mathbf{D}_3 \mathbf{D}_6,\\
  \mathbf{G}_{10} =&-\mathbf{D}_5^2+10 \mathbf{D}_0 \mathbf{D}_{10}-2 \mathbf{D}_{\infty}
                    \mathbf{D}_{10}+2 \mathbf{D}_6 \mathbf{D}_4+\\
                   &-2 \mathbf{D}_3 \mathbf{D}_7-4 \mathbf{D}_9 \mathbf{D}_1-2
                    \mathbf{D}_{12} \mathbf{D}_{11},\\
  \mathbf{G}_{11} =&-2 \mathbf{D}_{12}^2+6 \mathbf{D}_0 \mathbf{D}_{11}+4 \mathbf{D}_{\infty}
                    \mathbf{D}_{11}-2 \mathbf{D}_9 \mathbf{D}_2+\\
                   &-2 \mathbf{D}_5 \mathbf{D}_6-2 \mathbf{D}_7 \mathbf{D}_4-2
                    \mathbf{D}_3 \mathbf{D}_8,\\
  \mathbf{G}_{12} =&-\mathbf{D}_6^2+10 \mathbf{D}_0 \mathbf{D}_{12}-2 \mathbf{D}_{\infty}
                    \mathbf{D}_{12}+2 \mathbf{D}_2 \mathbf{D}_{10}+\\
                   &-2 \mathbf{D}_1 \mathbf{D}_{11}-4 \mathbf{D}_3 \mathbf{D}_9-2
                    \mathbf{D}_4 \mathbf{D}_8.
\endaligned\right.\eqno{(1.5)}$$
Here, the senary cubic forms (cubic forms in six variables)
$\mathbf{D}_j$ $(j=0$, $1$, $\ldots$, $12$, $\infty)$ are given as follows:
$$\left\{\aligned
  \mathbf{D}_0 &=z_1 z_2 z_3,\\
  \mathbf{D}_1 &=2 z_2 z_3^2+z_2^2 z_6-z_4^2 z_5+z_1 z_5 z_6,\\
  \mathbf{D}_2 &=-z_6^3+z_2^2 z_4-2 z_2 z_5^2+z_1 z_4 z_5+3 z_3 z_5 z_6,\\
  \mathbf{D}_3 &=2 z_1 z_2^2+z_1^2 z_5-z_4 z_6^2+z_3 z_4 z_5,\\
  \mathbf{D}_4 &=-z_2^2 z_3+z_1 z_6^2-2 z_4^2 z_6-z_1 z_3 z_5,\\
  \mathbf{D}_5 &=-z_4^3+z_3^2 z_5-2 z_3 z_6^2+z_2 z_5 z_6+3 z_1 z_4 z_6,\\
  \mathbf{D}_6 &=-z_5^3+z_1^2 z_6-2 z_1 z_4^2+z_3 z_4 z_6+3 z_2 z_4 z_5,\\
  \mathbf{D}_7 &=-z_2^3+z_3 z_4^2-z_1 z_3 z_6-3 z_1 z_2 z_5+2 z_1^2 z_4,\\
  \mathbf{D}_8 &=-z_1^3+z_2 z_6^2-z_2 z_3 z_5-3 z_1 z_3 z_4+2 z_3^2 z_6,\\
  \mathbf{D}_9 &=2 z_1^2 z_3+z_3^2 z_4-z_5^2 z_6+z_2 z_4 z_6,\\
  \mathbf{D}_{10} &=-z_1 z_3^2+z_2 z_4^2-2 z_4 z_5^2-z_1 z_2 z_6,\\
  \mathbf{D}_{11} &=-z_3^3+z_1 z_5^2-z_1 z_2 z_4-3 z_2 z_3 z_6+2 z_2^2 z_5,\\
  \mathbf{D}_{12} &=-z_1^2 z_2+z_3 z_5^2-2 z_5 z_6^2-z_2 z_3 z_4,\\
  \mathbf{D}_{\infty}&=z_4 z_5 z_6.
\endaligned\right.\eqno{(1.6)}$$
Then $\delta_{\infty}$, $\delta_{\nu}$ for $\nu=0, \ldots, 12$ are the
roots of a polynomial of degree fourteen. The corresponding equation is
not the Jacobian equation. Now, a family of invariants for $G$ is given
as follows: put
$$\Phi_4=\sum_{\nu=0}^{12} w_{\nu}+w_{\infty}, \quad
  \Phi_8=\sum_{\nu=0}^{12} w_{\nu}^2+w_{\infty}^2, \eqno{(1.7)}$$
$$\Phi_{12}=-\frac{1}{13 \cdot 52} \left(\sum_{\nu=0}^{12} \delta_{\nu}^2+
  \delta_{\infty}^2\right), \quad
  {\Phi}_{12}^{\prime}=-\frac{1}{13 \cdot 30} \left(\sum_{\nu=0}^{12}
  w_{\nu}^3+w_{\infty}^3\right), \eqno{(1.8)}$$
$$\Phi_{16}=\sum_{\nu=0}^{12} w_{\nu}^4+w_{\infty}^4, \quad
  \Phi_{18}=\frac{1}{13 \cdot 6} \left(\sum_{\nu=0}^{12} \delta_{\nu}^3
  +\delta_{\infty}^3\right), \eqno{(1.9)}$$
$$\Phi_{20}=\frac{1}{13 \cdot 25} \left(\sum_{\nu=0}^{12} w_{\nu}^5+
  w_{\infty}^5\right), \quad
  \Phi_{30}=-\frac{1}{13 \cdot 1315} \left(\sum_{\nu=0}^{12}
  \delta_{\nu}^5+\delta_{\infty}^5\right),\eqno{(1.10)}$$
and $x_i(z)=\eta(z) a_i(z)$ $(1 \leq i \leq 6)$, where
$$\left\{\aligned
  a_1(z) &:=e^{-\frac{11 \pi i}{26}} \theta \begin{bmatrix}
            \frac{11}{13}\\ 1 \end{bmatrix}(0, 13z),\\
  a_2(z) &:=e^{-\frac{7 \pi i}{26}} \theta \begin{bmatrix}
            \frac{7}{13}\\ 1 \end{bmatrix}(0, 13z),\\
  a_3(z) &:=e^{-\frac{5 \pi i}{26}} \theta \begin{bmatrix}
            \frac{5}{13}\\ 1 \end{bmatrix}(0, 13z),\\
  a_4(z) &:=-e^{-\frac{3 \pi i}{26}} \theta \begin{bmatrix}
            \frac{3}{13}\\ 1 \end{bmatrix}(0, 13z),\\
  a_5(z) &:=e^{-\frac{9 \pi i}{26}} \theta \begin{bmatrix}
            \frac{9}{13}\\ 1 \end{bmatrix}(0, 13z),\\
  a_6(z) &:=e^{-\frac{\pi i}{26}} \theta \begin{bmatrix}
            \frac{1}{13}\\ 1 \end{bmatrix}(0, 13z)
\endaligned\right.\eqno{(1.11)}$$
are theta constants of order $13$ and
$\eta(z):=q^{\frac{1}{24}} \prod_{n=1}^{\infty} (1-q^n)$ with
$q=e^{2 \pi i z}$ is the Dedekind eta function which are all
defined in the upper-half plane
$\mathbb{H}=\{ z \in \mathbb{C}: \text{Im}(z)>0 \}$.
Our main theorem is the following:

\textbf{Theorem 1.1.} {\it The $G$-invariant polynomials $\Phi_4$,
$\ldots$, $\Phi_{30}$ in $x_1(z)$, $\ldots$, $x_6(z)$ can be identified
with modular forms as follows$:$
$$\left\{\aligned
  \Phi_{4}(x_1(z), \ldots, x_6(z)) &=0,\\
  \Phi_{8}(x_1(z), \ldots, x_6(z)) &=0,\\
  \Phi_{12}(x_1(z), \ldots, x_6(z)) &=\Delta(z),\\
  \Phi_{12}^{\prime}(x_1(z), \ldots, x_6(z)) &=\Delta(z),\\
  \Phi_{16}(x_1(z), \ldots, x_6(z)) &=0,\\
  \Phi_{18}(x_1(z), \ldots, x_6(z)) &=\Delta(z) E_6(z),\\
  \Phi_{20}(x_1(z), \ldots, x_6(z)) &=\eta(z)^8 \Delta(z) E_4(z),\\
  \Phi_{30}(x_1(z), \ldots, x_6(z)) &=\Delta(z)^2 E_6(z).
\endaligned\right.\eqno{(1.12)}$$}

  Theorem 1.1 has many consequences. The first one comes from
the theory of singularities: there exists at least two kinds
of constructions of the equation of the $E_8$-singularity: one is
given by the icosahedral group in the celebrated book of Klein
\cite{K}, i.e., the icosahedral singularity (see \cite{Br3}, p. 107),
the other is given by the group $\text{PSL}(2, 13)$.

\textbf{Theorem 1.2} (A different construction of the $E_8$-singularity:
from $X(13)$ to $E_8$). {\it The equation of the $E_8$-singularity can
be constructed from the modular curve $X(13)$ as follows$:$
$$\Phi_{20}^3-\Phi_{30}^2=1728 \Phi_{12}^5, \eqno{(1.13)}$$
where $\Phi_j=\Phi_j(x_1(z), \ldots, x_6(z))$ for $j=12, 20$ and
$30$. As polynomials in six variables $z_1, \ldots, z_6$,
$\Phi_{12}$, $\Phi_{20}$ and $\Phi_{30}$ are $G$-invariant
polynomials.}

  In fact, in his talk at ICM 1970 \cite{Br2}, Brieskorn showed
how to construct the singularity of type $ADE$ directly from
the simple complex Lie group of the same type. At the end of
that paper \cite{Br2} Brieskorn says:``Thus we see that there
is a relation between exotic spheres, the icosahedron and $E_8$.
But I still do not understand why the regular polyhedra come in.''
(see also \cite{Gr}, \cite{GrP} and \cite{Br4}). As a consequence,
Theorem 1.2 shows that the $E_8$ root lattice is not necessarily
constructed from the icosahedron. That is, the icosahedron does
not necessarily appear in the triple (exotic spheres, icosahedron,
$E_8$) of Brieskorn \cite{Br2}. The group $\text{PSL}(2, 13)$ can
take its place and there is the other triple (exotic spheres,
$\text{PSL}(2, 13)$, $E_8$). The higher dimensional liftings of
these two distinct groups and modular interpretations on the
equation of the $E_8$-singularity give the same Milnor's standard
generator of $\Theta_7$.

  The second consequence of Theorem 1.1 comes from differential
topology of $4$-manifolds. The manifold $\mathbb{CP}^2 \# 8
\overline{\mathbb{CP}}^2$ has two distinct differentiable structures,
both of which come from algebraic surfaces: one is the eight-fold
blow-up of the projective plane carrying the standard smooth structure,
the other is the Barlow surface \cite{Ba}, which is a simply connected
minimal surface of general type with $q=p_g=0$ and $K^2=1$. In fact, up
to now, only two kinds of such surfaces are known: the first example is
the Barlow surface \cite{Ba}, the second examples are given by Lee and
Park in the appendix of \cite{LP}, both of them are simply connected,
minimal, complex surfaces of general type with $p_g=0$ and $K^2=1$. The
Barlow surface comes from a certain Hilbert modular surface associated to
the icosahedral group (see \cite{GZ}). On the other hand, the Lee-Park
surfaces are constructed by a rational blow-down surgery and a
$\mathbb{Q}$-Gorenstein smoothing theory (see \cite{LP}). The Barlow
surface is homeomorphic but not diffeomorphic to
$\mathbb{CP}^2 \# 8 \overline{\mathbb{CP}}^2$ (see \cite{Ko}). We
prove the following:

\textbf{Theorem 1.3} (A different construction of the Barlow surface:
from $X(13)$ to $E_8$). {\it The Barlow surface can be constructed
from the modular curve $X(13)$.}

  The third consequence of Theorem 1.1 comes from representation
theory and classical algebraic geometry: the automorphism group
of the configuration of $120$ tritangent planes of Bring's curve
is the quotient of the Weyl group $W(E_8)$ of the root system of
type $E_8$ by the normal subgroup $\{ \pm 1 \}$ with order
$2^{13} \cdot 3^5 \cdot 5^2 \cdot 7$, which is just the group
$G_{8, 2}$ studied by Coble in p. 356 of his paper \cite{Co}.
Note that Bring's curve has an analogue: Fricke's octavic curve,
both of them arise from the resolution of the equation of the
fifth degree (see \cite{K}). Hence, they are intimately connected
with the icosahedron. As a consequence of Theorem 1.1, we prove
the following:

\textbf{Theorem 1.4} (A different construction of Bring's curve
and Fricke's octavic curve: from $X(13)$ to $E_8$).  {\it Both
Bring's curve and Fricke's octavic curve can be constructed from
the modular curve $X(13)$.}

  Theorem 1.2, Theorem 1.3 and Theorem 1.4 show that there exist
two distinct constructions of the $E_8$ root lattice from its three
different appearances. Hence, the fact that the $E_8$ root lattice
is not necessarily constructed from the icosahedron, can be realized
not only from the theory of singularities and exotic spheres, but also
from differential topology of $4$-manifolds, representation theory and
classical algebraic geometry.

  The fourth consequence of Theorem 1.1 comes from the explicit
construction of modular curves, which is a classical problem
studied by Klein (see \cite{K1}, \cite{K2}, \cite{KF1} and
\cite{KF2}).

\textbf{Problem 1.5.} Let $p \geq 7$ be a prime number. Give an
explicit construction of the modular curve $X(p)$ of level $p$
from the invariant theory for $\mathrm{PSL}(2, p)$ using projective
algebraic geometry.

  For $p=7$, the modular curve $X(7)$ is given by the celebrated
Klein quartic curve (see \cite{K1})
$$x^3 y+y^3 z+z^3 x=0.$$
For $p=11$, the modular curve $X(11)$ leads to the study of
the Klein cubic threefold (see \cite{K2})
$$v^2 w+w^2 x+x^2 y+y^2 z+z^2 v=0.$$
Following Klein's method for the cubic threefold, Adler and Ramanan
(see \cite{AR}) studied Problem 1.5 when $p$ is a prime congruent
to $3$ modulo $8$ by some cubic hypersurface invariant under
$\mathrm{PSL}(2, p)$. However, their method can not be valid
for $p=13$. As a consequence of Theorem 1.1, we find an explicit
construction of the modular curve $X(13)$. Let
$$\phi_{12}(z_1, \ldots, z_6)=\Phi_{12}(z_1, \ldots, z_6)-
  \Phi_{12}^{\prime}(z_1, \ldots, z_6).\eqno{(1.14)}$$

\textbf{Theorem 1.6} (An explicit construction of the modular
curve $X(13)$). {\it There is a morphism
$$\Phi: X(13) \to C \subset \mathbb{CP}^5$$
with $\Phi(z)=(x_1(z), \ldots, x_6(z))$, where $C$ is an algebraic
curve given by a family of $G$-invariant equations
$$\left\{\aligned
  \Phi_{4}(z_1, \ldots, z_6) &=0,\\
  \Phi_{8}(z_1, \ldots, z_6) &=0,\\
  \phi_{12}(z_1, \ldots, z_6) &=0,\\
  \Phi_{16}(z_1, \ldots, z_6) &=0,
\endaligned\right.\eqno{(1.15)}$$
where $\phi_{12}$ is given as in (1.14).}

  This paper consists of five sections. In section two, we revisit the
standard construction of the $E_8$ root lattice by means of the icosahedron.
This includes the $E_8$-singularity, the Barlow surface, Bring's curve
and Fricke's octavic curve. In section three, we explain the invariant
theory for $\text{PSL}(2, 13)$. In particular, we construct
the senary quadratic forms $\mathbf{A}_j$ ($0 \leq j \leq 6$), the
senary cubic forms $\mathbf{D}_j$ ($j=0, 1, \ldots, 12, \infty$) and
the senary sextic forms $\mathbf{G}_j$ ($0 \leq j \leq 12$). From
$\mathbf{A}_j$ we construct the Jacobian equation of degree fourteen. From
$\mathbf{D}_j$ and $\mathbf{G}_j$ we construct another equation of degree
fourteen. Combining Jacobian equation with that equation, we obtain a
family of polynomials which are invariant under the action of
$\text{PSL}(2, 13)$. Together with theta constants of order thirteen,
this gives the modular parametrization of these invariant polynomials.
Therefore, we obtain Theorem 1.1. In section four, we give three
constructions of the $E_8$ root lattice, and prove Theorem 1.2, Theorem
1.3 and Theorem 1.4. In section five, we give an explicit construction
of the modular curve $X(13)$.

\vskip 0.3 cm

\textbf{Acknowledgements}. The author would like to thank Jean-Pierre
Serre for his very detailed and helpful comments and his patience. The
author also thanks Pierre Deligne for his helpful comments.

\begin{center}
{\large\bf 2. Standard construction: from the icosahedron to $E_8$}
\end{center}

\noindent{\bf 2.1. $E_8$-singularity: from the icosahedron to $E_8$}

  Let us recall some classical result on the relation between the
icosahedron and the $E_8$-singularity (see \cite{Mc}). Starting with
the polynomial invariants of the finite subgroup of $\text{SL}(2, \mathbb{C})$,
a surface is defined from the single syzygy which
relates the three polynomials in two variables. This surface has a singularity
at the origin; the singularity can be resolved by constructing a smooth surface
which is isomorphic to the original one except for a set of component curves
which form the pre-image of the origin. The components form a Dynkin curve and
the matrix of their intersections is the negative of the Cartan matrix for the
appropriate Lie algebra. The Dynkin curve is the dual of the Dynkin graph. For
example, if $\Gamma$ is the binary icosahedral group, the corresponding Dynkin
curve is that of $E_8$, and $\mathbb{C}^2/\Gamma \subset \mathbb{C}^3$ is the set
of zeros of the equation
$$x^2+y^3+z^5=0.\eqno{(2.1)}$$
The link of this $E_8$-singularity, the Poincar\'{e} homology $3$-sphere (see
\cite{KS}), has a higher dimensional lifting:
$$z_1^5+z_2^3+z_3^2+z_4^2+z_5^2=0, \quad \sum_{i=1}^{5} z_i \overline{z_i}=1,
  \quad z_i \in \mathbb{C} \quad (1 \leq i \leq 5),\eqno{(2.2)}$$
which is the Brieskorn description of one of Milnor's exotic $7$-dimensional
spheres. In fact, it is an exotic $7$-sphere representing Milnor's standard
generator of $\Theta_7$ (see \cite{Br1}, \cite{Br2} and \cite{Hi}).

  In his celebrated book \cite{K}, Klein gave a parametric solution of the
above singularity (2.1) by homogeneous polynomials $T$, $H$, $f$ in two
variables of degrees $30$, $20$, $12$ with integral coefficients, where
$$f=z_1 z_2 (z_1^{10}+11 z_1^5 z_2^5-z_2^{10}),$$
$$H=\frac{1}{121} \begin{vmatrix}
    \frac{\partial^2 f}{\partial z_1^2} &
    \frac{\partial^2 f}{\partial z_1 \partial z_2}\\
    \frac{\partial^2 f}{\partial z_2 \partial z_1} &
    \frac{\partial^2 f}{\partial z_2^2}
    \end{vmatrix}
  =-(z_1^{20}+z_2^{20})+228 (z_1^{15} z_2^5-z_1^5 z_2^{15})
   -494 z_1^{10} z_2^{10},$$
$$T=-\frac{1}{20} \begin{vmatrix}
    \frac{\partial f}{\partial z_1} &
    \frac{\partial f}{\partial z_2}\\
    \frac{\partial H}{\partial z_1} &
    \frac{\partial H}{\partial z_2}
    \end{vmatrix}
  =(z_1^{30}+z_2^{30})+522 (z_1^{25} z_2^5-z_1^5 z_2^{25})
   -10005 (z_1^{20} z_2^{10}+z_1^{10} z_2^{20}).$$
They satisfy the famous (binary) icosahedral equation
$$T^2+H^3=1728 f^5.\eqno{(2.3)}$$
In fact, $f$, $H$ and $T$ are invariant polynomials under the action
of the binary icosahedral group. The above equation (2.3) is closely
related to Hermite's celebrated work (see \cite{He}, pp.5-12) on the
resolution of the quintic equations. Essentially the same relation had been
found a few years earlier by Schwarz (see \cite{Sch}), who considered three
polynomials $\varphi_{12}$, $\varphi_{20}$ and $\varphi_{30}$ whose roots
correspond to the vertices, the midpoints of the faces and the midpoints
of the edges of an icosahedron inscribed in the Riemann sphere. He obtained
the identity $\varphi_{20}^3-1728 \varphi_{12}^5=\varphi_{30}^2$.
Thus we see that from the very beginning there was a close relation between
the $E_8$-singularity and the icosahedron. Moreover, the icosahedral equation
(2.3) can be interpreted in terms of modular forms which was also known by
Klein (see \cite{KF1}, p. 631). Let $x_1(z)=\eta(z) a(z)$ and
$x_2(z)=\eta(z) b(z)$, where
$$a(z)=e^{-\frac{3 \pi i}{10}} \theta \begin{bmatrix}
       \frac{3}{5}\\ 1 \end{bmatrix}(0, 5z), \quad
  b(z)=e^{-\frac{\pi i}{10}} \theta \begin{bmatrix}
       \frac{1}{5}\\ 1 \end{bmatrix}(0, 5z)$$
are theta constants of order five and
$\eta(z):=q^{\frac{1}{24}} \prod_{n=1}^{\infty} (1-q^n)$ with
$q=e^{2 \pi i z}$ is the Dedekind eta function which are all
defined in the upper-half plane
$\mathbb{H}=\{ z \in \mathbb{C}: \text{Im}(z)>0 \}$. Then
$$\left\{\aligned
  f(x_1(z), x_2(z)) &=-\Delta(z),\\
  H(x_1(z), x_2(z)) &=-\eta(z)^8 \Delta(z)E_4(z),\\
  T(x_1(z), x_2(z)) &=\Delta(z)^2 E_6(z),
\endaligned\right.$$
where
$$E_4(z):=\frac{1}{2} \sum_{m, n \in \mathbb{Z}, (m, n)=1}
          \frac{1}{(mz+n)^4}, \quad
  E_6(z):=\frac{1}{2} \sum_{m, n \in \mathbb{Z}, (m, n)=1}
          \frac{1}{(mz+n)^6}$$
are Eisenstein series of weight $4$ and $6$, and
$\Delta(z)=\eta(z)^{24}$ is the discriminant. The relations
$$j(z):=\frac{E_4(z)^3}{\Delta(z)}=\frac{H(x_1(z), x_2(z))^3}
        {f(x_1(z), x_2(z))^5},$$
$$j(z)-1728=\frac{E_6(z)^2}{\Delta(z)}=-\frac{T(x_1(z), x_2(z))^2}
            {f(x_1(z), x_2(z))^5}$$
give the icosahedral equation (2.3) in terms of theta constants of
order five.

\vskip 0.3 cm

\noindent{\bf 2.2. The Barlow surface: from the icosahedron to $E_8$}

  In \cite{Ko}, Kotschick showed that the manifold $\mathbb{CP}^2 \#
8 \overline{\mathbb{CP}}^2$ has two distinct differentiable structures,
both of which come from algebraic surfaces. The surfaces are the eight-fold
blow-up of the projective plane carrying the standard smooth structure and
the Barlow surface (see \cite{Ba}), which is a simply-connected minimal
surface of general type with $q=p_g=0$ and $K^2=1$. Kotschick proved that
the Barlow surface $B$ is homeomorphic but not diffeomorphic to
$\mathbb{CP}^2 \# 8 \overline{\mathbb{CP}}^2$.

  The Barlow surface $B$ (see \cite{Ba}) is obtained as the minimal
desingularization of $Y/D_{10}$, where $Y$ is a certain Hilbert modular
surface (see \cite{GZ}) and $D_{10}$ acts with finite fixed locus. Let
$\mathbb{H}=\{ z \in \mathbb{C}: \mathrm{Im}(z)>0 \}$ be the upper half-plane.
Denote by $K$ the real quadratic field $\mathbb{Q}(\sqrt{21})$ and by
$\mathcal{O}_K$ its ring of integers. The Hilbert modular group
$\mathrm{SL}(2, \mathcal{O}_K)$ acts on $\mathbb{H} \times \mathbb{H}$
by
$$\begin{pmatrix}
   \alpha & \beta\\
   \gamma & \delta
  \end{pmatrix} \circ (z_1, z_2)
 =\left(\frac{\alpha z_1+\beta}{\gamma z_1+\delta},
        \frac{\alpha^{\prime} z_2+\beta^{\prime}}
        {\gamma^{\prime} z_2+\delta^{\prime}}\right),$$
where $x \mapsto x^{\prime}$ denotes conjugation over $\mathbb{Q}$
in $K$. Consider the $2$-congruence subgroup $\Gamma \subset
\mathrm{SL}(2, \mathcal{O}_K)$, where
$$\Gamma=\left\{ \begin{pmatrix}
          \alpha & \beta\\
          \gamma & \delta
          \end{pmatrix} \in \mathrm{SL}(2, \mathcal{O}_K):
          \alpha \equiv \delta \equiv 1,
          \beta \equiv \gamma \equiv 0 (\mbox{mod $2$}) \right\}$$
is the principal congruence subgroup of $\mathrm{SL}(2, \mathcal{O}_K)$
for the prime ideal generated by $2$. The surface $(\mathbb{H} \times
\mathbb{H})/\mathrm{SL}(2, \mathcal{O}_K)$ is the quotient of
$(\mathbb{H} \times \mathbb{H})/\Gamma$ by the group
$$\mathrm{SL}(2, \mathcal{O}_K)/\Gamma \cong
  \mathrm{SL}(2, \mathcal{O}_K/2 \mathcal{O}_K) \cong
  \mathrm{SL}(2, \mathbb{F}_4) \cong A_5.\eqno{(2.4)}$$
The surface $Y$ is the minimal desingularisation of the resolution of
the compactification of $(\mathbb{H} \times \mathbb{H})/\Gamma$. It
is a simply-connected surface of general type with $p_g=4$ and
$K^2=10$, It can be proved (see \cite{GZ}) that the canonical map of
$Y$ is $2:1$ onto the $20$-nodal quintic $Q \subset \mathbb{CP}^4$
given by
$$\sum_{i=0}^{4} z_i=0, \quad
  \sum_{i=0}^{4} z_i^5-\frac{5}{4} \sum_{i=0}^{4} z_i^2 \cdot
  \sum_{i=0}^{4} z_i^3=0.\eqno{(2.5)}$$
The icosahedral group $A_5$ acts on $\mathbb{CP}^4$ by the standard
action on the coordinates. The quintic $Q$ is $A_5$-invariant and
its $20$ nodes are the $A_5$-orbit of the point $(2, 2, 2,
-3-\sqrt{-7}, -3+\sqrt{-7})$. The $A_5$-action on $Q$ is covered
by an action on $Y$, so that we have an action of $A_5 \times
\mathbb{Z}/2 \mathbb{Z}=A_5 \cup A_5 \sigma$ on $Y$, where the
generator $\sigma \in \mathbb{Z}/2 \mathbb{Z}$ is the covering
involution. Elements of $A_5 \times \mathbb{Z}/2 \mathbb{Z}$
acting on $Y$ are denoted like the corresponding elements acting
on $Q$. Let $\Phi: Y \to Q$ be the quotient map. The main result
from \cite{Ba} can be summarized as follows:

\textbf{Proposition 2.1} (see \cite{Ba} and \cite{Ko}). {\it Let
$\alpha=(02)(34) \sigma$, $\beta=(01234)$. Then $\beta$ acts freely
on $Y$ and $\alpha$ has $4$ fixed points. The resolution of the nodes
of $Y/D_{10}$, where $D_{10}=\langle \alpha, \beta \rangle$, gives a
minimal surface $B$ of general type with $\pi_1=0$, $q=p_g=0$ and
$K^2=1$.}

\vskip 0.3 cm

\noindent{\bf 2.3. Bring's curve and Fricke's octavic curve: from the
                icosahedron to $E_8$}

  Both Bring's curve and Fricke's octavic curve arise from the
resolution of the equation of the fifth degree (see \cite{K}).
Hence, they are intimately connected with the icosahedron.

  Let us recall some basic facts from classical enumerative geometry
(see \cite{HM}). We will use the notation that $[n|a_1, a_2, \ldots,
a_k]$ means the (not necessarily complete) intersection of $k$
polynomials of degrees $a_1, \ldots, a_k$ respectively in
$\mathbb{P}^n$. In 1863, Clebsch found that the canonical sextic
curve of genus four has exactly $120$ tritangent planes (i.e., planes
which are tangent to the curve at precisely three points). This curve
can be realized as $[4|1, 2, 3]$, i.e., the intersection of a hyperplane,
a quadric and a cubic in Fermat form in the five homogeneous
coordinates of $\mathbb{P}^4$, giving us the so-called Bring's curve
(see \cite{Ed1} and \cite{Ed2}).

  Specially, Bring's curve can be realized as the Fermat cubic,
sliced by the Fermat quadric, and then the line, in the homogeneous
coordinates of $\mathbb{P}^4$:
$$\mathcal{B}=\left\{ \sum_{i=0}^{4} x_i^3=\sum_{i=0}^{4} x_i^2
             =\sum_{i=0}^{4} x_i=0 \right\} \subset \mathbb{P}^4.\eqno{(2.6)}$$
This classic result is well-understood in terms of del Pezzo surfaces of degree
one (see \cite{M}). The canonical model of a del Pezzo surface of degree one
is the double cover of a quadratic cone, branched over a canonical space curve
of genus $4$ and degree $6$ given by the complete intersection of the cone with
a unique cubic surface. The $240$ lines on the del Pezzo surface arise in pairs
from the $120$ tritangent planes to the canonical curve, which can be identified
with its odd theta characteristics.

  The automorphism group of the $240$ lines is the Weyl group $W(E_8)$ of the
root system of type $E_8$. Its order is $2^{14} \cdot 3^5 \cdot 5^2 \cdot 7$.
The automorphism group of the $120$ tritangent planes is the quotient by the
normal subgroup $\{ \pm 1 \}$ with order $2^{13} \cdot 3^5 \cdot 5^2 \cdot 7$,
which is just the group $G_{8, 2}$ studied by Coble in p. 356 of his paper
\cite{Co}. In fact, Coble proved that $G_{8, 2}$ is isomorphic with the group
of the tritangent planes of a space sextic of genus $4$ on a quadric cone (see
\cite{Co}, p. 359).

  Note that Bring's curve has a natural modular interpretation
(see \cite{Do}, p. 500). Namely, it is isomorphic to the modular
curve $\overline{\mathbb{H}/\Gamma}$, where $\Gamma=\Gamma_0(2)
\cap \Gamma(5)$. It is also realized as the curve of fixed points
of the Bertini involution on the del Pezzo surface of degree one
obtained from the elliptic modular surface $S(5)$ of level $5$ by
blowing down the zero section.

  Moreover, there is even a correspondence between the above
classical enumerative geometrical problem and the Monster simple
group $\mathbb{M}$ due to \cite{HM}, observation 1: for the
Monster $\mathbb{M}$, we have the following sums for the cusp
numbers $C_g$ over the $172$ rational conjugacy classes:
$$\sum_{g} C_g=360=3 \cdot 120, \quad
  \sum_{g} C_g^2=1024=2^{10}.$$
The $360$ is thrice $120$, which is the number of tritangent
planes to Bring's curve. Furthermore, in analogy to Bring's
sextic curve, there is the octavic of Fricke of genus nine (see
\cite{Ed3} and \cite{F}), the Fermat $[4|1, 2, 4]$ defined as
$$\mathcal{F}=\left\{ \sum_{i=0}^{4} x_i^4=\sum_{i=0}^{4} x_i^2
             =\sum_{i=0}^{4} x_i=0 \right\} \subset \mathbb{P}^4.\eqno{(2.7)}$$
The number of tritangent planes on $\mathcal{F}$ is precisely
$2048=2 \cdot 1024$, twice the sum of square of the cusps.

\begin{center}
{\large\bf 3. Modular curve $X(13)$ and invariant theory for $\text{PSL}(2, 13)$}
\end{center}

  At first, we will study the six-dimensional representation of the
finite simple group $\text{PSL}(2, 13)$ of order $1092$, which acts
on the five-dimensional projective space
$\mathbb{P}^5=\{ (z_1, z_2, z_3, z_4, z_5, z_6): z_i \in \mathbb{C}
 \quad (i=1, 2, 3, 4, 5, 6) \}$. This representation is defined over
the cyclotomic field $\mathbb{Q}(e^{\frac{2 \pi i}{13}})$. Put
$$S=-\frac{1}{\sqrt{13}} \begin{pmatrix}
  \zeta^{12}-\zeta & \zeta^{10}-\zeta^3 & \zeta^4-\zeta^9 &
  \zeta^5-\zeta^8 & \zeta^2-\zeta^{11} & \zeta^6-\zeta^7\\
  \zeta^{10}-\zeta^3 & \zeta^4-\zeta^9 & \zeta^{12}-\zeta &
  \zeta^2-\zeta^{11} & \zeta^6-\zeta^7 & \zeta^5-\zeta^8\\
  \zeta^4-\zeta^9 & \zeta^{12}-\zeta & \zeta^{10}-\zeta^3 &
  \zeta^6-\zeta^7 & \zeta^5-\zeta^8 & \zeta^2-\zeta^{11}\\
  \zeta^5-\zeta^8 & \zeta^2-\zeta^{11} & \zeta^6-\zeta^7 &
  \zeta-\zeta^{12} & \zeta^3-\zeta^{10} & \zeta^9-\zeta^4\\
  \zeta^2-\zeta^{11} & \zeta^6-\zeta^7 & \zeta^5-\zeta^8 &
  \zeta^3-\zeta^{10} & \zeta^9-\zeta^4 & \zeta-\zeta^{12}\\
  \zeta^6-\zeta^7 & \zeta^5-\zeta^8 & \zeta^2-\zeta^{11} &
  \zeta^9-\zeta^4 & \zeta-\zeta^{12} & \zeta^3-\zeta^{10}
\end{pmatrix}\eqno{(3.1)}$$
and
$$T=\text{diag}(\zeta^7, \zeta^{11}, \zeta^8, \zeta^6,
    \zeta^2, \zeta^5),\eqno{(3.2)}$$
where $\zeta=\exp(2 \pi i/13)$. We have
$$S^2=T^{13}=(ST)^3=1.\eqno{(3.3)}$$
Let $G=\langle S, T \rangle$, then $G \cong \text{PSL}(2, 13)$
(see \cite{Y1}, Theorem 3.1).

  Put $\theta_1=\zeta+\zeta^3+\zeta^9$,
$\theta_2=\zeta^2+\zeta^6+\zeta^5$,
$\theta_3=\zeta^4+\zeta^{12}+\zeta^{10}$,
and $\theta_4=\zeta^8+\zeta^{11}+\zeta^7$. We find that
$$\left\{\aligned
  &\theta_1+\theta_2+\theta_3+\theta_4=-1,\\
  &\theta_1 \theta_2+\theta_1 \theta_3+\theta_1 \theta_4+
   \theta_2 \theta_3+\theta_2 \theta_4+\theta_3 \theta_4=2,\\
  &\theta_1 \theta_2 \theta_3+\theta_1 \theta_2 \theta_4+
   \theta_1 \theta_3 \theta_4+\theta_2 \theta_3 \theta_4=4,\\
  &\theta_1 \theta_2 \theta_3 \theta_4=3.
\endaligned\right.$$
Hence, $\theta_1$, $\theta_2$, $\theta_3$ and $\theta_4$ satisfy
the quartic equation $z^4+z^3+2 z^2-4z+3=0$,
which can be decomposed as two quadratic equations
$$\left(z^2+\frac{1+\sqrt{13}}{2} z+\frac{5+\sqrt{13}}{2}\right)
  \left(z^2+\frac{1-\sqrt{13}}{2} z+\frac{5-\sqrt{13}}{2}\right)=0$$
over the real quadratic field $\mathbb{Q}(\sqrt{13})$. Therefore, the
four roots are given as follows:
$$\left\{\aligned
  \theta_1=\frac{1}{4} \left(-1+\sqrt{13}+\sqrt{-26+6 \sqrt{13}}\right),\\
  \theta_2=\frac{1}{4} \left(-1-\sqrt{13}+\sqrt{-26-6 \sqrt{13}}\right),\\
  \theta_3=\frac{1}{4} \left(-1+\sqrt{13}-\sqrt{-26+6 \sqrt{13}}\right),\\
  \theta_4=\frac{1}{4} \left(-1-\sqrt{13}-\sqrt{-26-6 \sqrt{13}}\right).
\endaligned\right.$$
Moreover, we find that
$$\left\{\aligned
  \theta_1+\theta_3+\theta_2+\theta_4 &=-1,\\
  \theta_1+\theta_3-\theta_2-\theta_4 &=\sqrt{13},\\
  \theta_1-\theta_3-\theta_2+\theta_4 &=-\sqrt{-13+2 \sqrt{13}},\\
  \theta_1-\theta_3+\theta_2-\theta_4 &=\sqrt{-13-2 \sqrt{13}}.
\endaligned\right.$$

  Let us study the action of $S T^{\nu}$ on $\mathbb{P}^5$, where
$\nu=0, 1, \ldots, 12$. Put
$$\alpha=\zeta+\zeta^{12}-\zeta^5-\zeta^8, \quad
   \beta=\zeta^3+\zeta^{10}-\zeta^2-\zeta^{11}, \quad
   \gamma=\zeta^9+\zeta^4-\zeta^6-\zeta^7.$$
We find that
$$\aligned
  &13 ST^{\nu}(z_1) \cdot ST^{\nu}(z_4)\\
=&\beta z_1 z_4+\gamma z_2 z_5+\alpha z_3 z_6+\\
 &+\gamma \zeta^{\nu} z_1^2+\alpha \zeta^{9 \nu} z_2^2+\beta \zeta^{3 \nu} z_3^2
  -\gamma \zeta^{12 \nu} z_4^2-\alpha \zeta^{4 \nu} z_5^2-\beta \zeta^{10 \nu} z_6^2+\\
 &+(\alpha-\beta) \zeta^{5 \nu} z_1 z_2+(\beta-\gamma) \zeta^{6 \nu} z_2 z_3
  +(\gamma-\alpha) \zeta^{2 \nu} z_1 z_3+\\
 &+(\beta-\alpha) \zeta^{8 \nu} z_4 z_5+(\gamma-\beta) \zeta^{7 \nu} z_5 z_6
  +(\alpha-\gamma) \zeta^{11 \nu} z_4 z_6+\\
 &-(\alpha+\beta) \zeta^{\nu} z_3 z_4-(\beta+\gamma) \zeta^{9 \nu} z_1 z_5
  -(\gamma+\alpha) \zeta^{3 \nu} z_2 z_6+\\
 &-(\alpha+\beta) \zeta^{12 \nu} z_1 z_6-(\beta+\gamma) \zeta^{4 \nu} z_2 z_4
  -(\gamma+\alpha) \zeta^{10 \nu} z_3 z_5.
\endaligned$$
$$\aligned
  &13 ST^{\nu}(z_2) \cdot ST^{\nu}(z_5)\\
=&\gamma z_1 z_4+\alpha z_2 z_5+\beta z_3 z_6+\\
 &+\alpha \zeta^{\nu} z_1^2+\beta \zeta^{9 \nu} z_2^2+\gamma \zeta^{3 \nu} z_3^2
  -\alpha \zeta^{12 \nu} z_4^2-\beta \zeta^{4 \nu} z_5^2-\gamma \zeta^{10 \nu} z_6^2+\\
 &+(\beta-\gamma) \zeta^{5 \nu} z_1 z_2+(\gamma-\alpha) \zeta^{6 \nu} z_2 z_3
  +(\alpha-\beta) \zeta^{2 \nu} z_1 z_3+\\
 &+(\gamma-\beta) \zeta^{8 \nu} z_4 z_5+(\alpha-\gamma) \zeta^{7 \nu} z_5 z_6
  +(\beta-\alpha) \zeta^{11 \nu} z_4 z_6+\\
 &-(\beta+\gamma) \zeta^{\nu} z_3 z_4-(\gamma+\alpha) \zeta^{9 \nu} z_1 z_5
  -(\alpha+\beta) \zeta^{3 \nu} z_2 z_6+\\
 &-(\beta+\gamma) \zeta^{12 \nu} z_1 z_6-(\gamma+\alpha) \zeta^{4 \nu} z_2 z_4
  -(\alpha+\beta) \zeta^{10 \nu} z_3 z_5.
\endaligned$$
$$\aligned
  &13 ST^{\nu}(z_3) \cdot ST^{\nu}(z_6)\\
=&\alpha z_1 z_4+\beta z_2 z_5+\gamma z_3 z_6+\\
 &+\beta \zeta^{\nu} z_1^2+\gamma \zeta^{9 \nu} z_2^2+\alpha \zeta^{3 \nu} z_3^2
  -\beta \zeta^{12 \nu} z_4^2-\gamma \zeta^{4 \nu} z_5^2-\alpha \zeta^{10 \nu} z_6^2+\\
 &+(\gamma-\alpha) \zeta^{5 \nu} z_1 z_2+(\alpha-\beta) \zeta^{6 \nu} z_2 z_3
  +(\beta-\gamma) \zeta^{2 \nu} z_1 z_3+\\
 &+(\alpha-\gamma) \zeta^{8 \nu} z_4 z_5+(\beta-\alpha) \zeta^{7 \nu} z_5 z_6
  +(\gamma-\beta) \zeta^{11 \nu} z_4 z_6+\\
 &-(\gamma+\alpha) \zeta^{\nu} z_3 z_4-(\alpha+\beta) \zeta^{9 \nu} z_1 z_5
  -(\beta+\gamma) \zeta^{3 \nu} z_2 z_6+\\
 &-(\gamma+\alpha) \zeta^{12 \nu} z_1 z_6-(\alpha+\beta) \zeta^{4 \nu} z_2 z_4
  -(\beta+\gamma) \zeta^{10 \nu} z_3 z_5.
\endaligned$$
Note that $\alpha+\beta+\gamma=\sqrt{13}$, we find that
$$\aligned
  &\sqrt{13} \left[ST^{\nu}(z_1) \cdot ST^{\nu}(z_4)+
   ST^{\nu}(z_2) \cdot ST^{\nu}(z_5)+ST^{\nu}(z_3) \cdot ST^{\nu}(z_6)\right]\\
 =&(z_1 z_4+z_2 z_5+z_3 z_6)+(\zeta^{\nu} z_1^2+\zeta^{9 \nu} z_2^2+\zeta^{3 \nu} z_3^2)
  -(\zeta^{12 \nu} z_4^2+\zeta^{4 \nu} z_5^2+\zeta^{10 \nu} z_6^2)+\\
  &-2(\zeta^{\nu} z_3 z_4+\zeta^{9 \nu} z_1 z_5+\zeta^{3 \nu} z_2 z_6)
   -2(\zeta^{12 \nu} z_1 z_6+\zeta^{4 \nu} z_2 z_4+\zeta^{10 \nu} z_3 z_5).
\endaligned$$
Let
$$\varphi_{\infty}(z_1, z_2, z_3, z_4, z_5, z_6)=\sqrt{13} (z_1 z_4+z_2 z_5+z_3 z_6)
  \eqno{(3.4)}$$
and
$$\varphi_{\nu}(z_1, z_2, z_3, z_4, z_5, z_6)=\varphi_{\infty}(ST^{\nu}(z_1, z_2,
                                              z_3, z_4, z_5, z_6))\eqno{(3.5)}$$
for $\nu=0, 1, \ldots, 12$. Then
$$\aligned
  \varphi_{\nu}
=&(z_1 z_4+z_2 z_5+z_3 z_6)+\zeta^{\nu} (z_1^2-2 z_3 z_4)+\zeta^{4 \nu} (-z_5^2-2 z_2 z_4)+\\
 &+\zeta^{9 \nu} (z_2^2-2 z_1 z_5)+\zeta^{3 \nu} (z_3^2-2 z_2 z_6)+
   \zeta^{12 \nu} (-z_4^2-2 z_1 z_6)+\\
 &+\zeta^{10 \nu} (-z_6^2-2 z_3 z_5).
\endaligned\eqno{(3.6)}$$
This leads us to define the following senary quadratic forms
(quadratic forms in six variables):
$$\left\{\aligned
  \mathbf{A}_0 &=z_1 z_4+z_2 z_5+z_3 z_6,\\
  \mathbf{A}_1 &=z_1^2-2 z_3 z_4,\\
  \mathbf{A}_2 &=-z_5^2-2 z_2 z_4,\\
  \mathbf{A}_3 &=z_2^2-2 z_1 z_5,\\
  \mathbf{A}_4 &=z_3^2-2 z_2 z_6,\\
  \mathbf{A}_5 &=-z_4^2-2 z_1 z_6,\\
  \mathbf{A}_6 &=-z_6^2-2 z_3 z_5.
\endaligned\right.\eqno{(3.7)}$$
Hence,
$$\sqrt{13} ST^{\nu}(\mathbf{A}_0)=\mathbf{A}_0+\zeta^{\nu} \mathbf{A}_1
  +\zeta^{4 \nu} \mathbf{A}_2+\zeta^{9 \nu} \mathbf{A}_3+\zeta^{3 \nu}
  \mathbf{A}_4+\zeta^{12 \nu} \mathbf{A}_5+\zeta^{10 \nu} \mathbf{A}_6.
  \eqno{(3.8)}$$
Let $H:=Q^5 P^2 \cdot P^2 Q^6 P^8 \cdot Q^5 P^2 \cdot P^3 Q$ where
$P=S T^{-1} S$ and $Q=S T^3$. Then (see \cite{Y2}, p.27)
$$H=\begin{pmatrix}
  0 &  0 &  0 & 0 & 0 & 1\\
  0 &  0 &  0 & 1 & 0 & 0\\
  0 &  0 &  0 & 0 & 1 & 0\\
  0 &  0 & -1 & 0 & 0 & 0\\
 -1 &  0 &  0 & 0 & 0 & 0\\
  0 & -1 &  0 & 0 & 0 & 0
\end{pmatrix}.\eqno{(3.9)}$$
Note that $H^6=1$ and $H^{-1} T H=-T^4$. Thus,
$\langle H, T \rangle \cong \mathbb{Z}_{13} \rtimes \mathbb{Z}_6$.
Hence, it is a maximal subgroup of order $78$ of $G$ with index
$14$. We find that $\varphi_{\infty}^2$ is invariant under the
action of the maximal subgroup $\langle H, T \rangle$. Note that
$$\varphi_{\infty}=\sqrt{13} \mathbf{A}_0, \quad
  \varphi_{\nu}=\mathbf{A}_0+\zeta^{\nu} \mathbf{A}_1+
  \zeta^{4 \nu} \mathbf{A}_2+\zeta^{9 \nu} \mathbf{A}_3+
  \zeta^{3 \nu} \mathbf{A}_4+\zeta^{12 \nu} \mathbf{A}_5+
  \zeta^{10 \nu} \mathbf{A}_6$$
for $\nu=0, 1, \ldots, 12$. Let $w=\varphi^2$,
$w_{\infty}=\varphi_{\infty}^2$ and $w_{\nu}=\varphi_{\nu}^2$.
Then $w_{\infty}$, $w_{\nu}$ for $\nu=0, \ldots, 12$ form an
algebraic equation of degree fourteen, which is just the Jacobian
equation of degree fourteen (see \cite{K}, pp.161-162), whose roots
are these $w_{\nu}$ and $w_{\infty}$:
$$w^{14}+a_1 w^{13}+\cdots+a_{13} w+a_{14}=0.$$

  On the other hand, we have
$$\aligned
  &-13 \sqrt{13} ST^{\nu}(z_1) \cdot ST^{\nu}(z_2) \cdot ST^{\nu}(z_3)\\
 =&-r_4 (\zeta^{8 \nu} z_1^3+\zeta^{7 \nu} z_2^3+\zeta^{11 \nu} z_3^3)
   -r_2 (\zeta^{5 \nu} z_4^3+\zeta^{6 \nu} z_5^3+\zeta^{2 \nu} z_6^3)\\
  &-r_3 (\zeta^{12 \nu} z_1^2 z_2+\zeta^{4 \nu} z_2^2 z_3+\zeta^{10 \nu} z_3^2 z_1)
   -r_1 (\zeta^{\nu} z_4^2 z_5+\zeta^{9 \nu} z_5^2 z_6+\zeta^{3 \nu} z_6^2 z_4)\\
  &+2 r_1 (\zeta^{3 \nu} z_1 z_2^2+\zeta^{\nu} z_2 z_3^2+\zeta^{9 \nu} z_3 z_1^2)
   -2 r_3 (\zeta^{10 \nu} z_4 z_5^2+\zeta^{12 \nu} z_5 z_6^2+\zeta^{4 \nu} z_6 z_4^2)\\
  &+2 r_4 (\zeta^{7 \nu} z_1^2 z_4+\zeta^{11 \nu} z_2^2 z_5+\zeta^{8 \nu} z_3^2 z_6)
   -2 r_2 (\zeta^{6 \nu} z_1 z_4^2+\zeta^{2 \nu} z_2 z_5^2+\zeta^{5 \nu} z_3 z_6^2)+\\
  &+r_1 (\zeta^{3 \nu} z_1^2 z_5+\zeta^{\nu} z_2^2 z_6+\zeta^{9 \nu} z_3^2 z_4)
   +r_3 (\zeta^{10 \nu} z_2 z_4^2+\zeta^{12 \nu} z_3 z_5^2+\zeta^{4 \nu} z_1 z_6^2)+\\
  &+r_2 (\zeta^{6 \nu} z_1^2 z_6+\zeta^{2 \nu} z_2^2 z_4+\zeta^{5 \nu} z_3^2 z_5)
   +r_4 (\zeta^{7 \nu} z_3 z_4^2+\zeta^{11 \nu} z_1 z_5^2+\zeta^{8 \nu} z_2 z_6^2)+\\
  &+r_0 z_1 z_2 z_3+r_{\infty} z_4 z_5 z_6+\\
  &-r_4 (\zeta^{11 \nu} z_1 z_2 z_4+\zeta^{8 \nu} z_2 z_3 z_5+\zeta^{7 \nu} z_1 z_3 z_6)+\\
  &+r_2 (\zeta^{2 \nu} z_1 z_4 z_5+\zeta^{5 \nu} z_2 z_5 z_6+\zeta^{6 \nu} z_3 z_4 z_6)+\\
  &-3 r_4 (\zeta^{7 \nu} z_1 z_2 z_5+\zeta^{11 \nu} z_2 z_3 z_6+\zeta^{8 \nu} z_1 z_3 z_4)+\\
  &+3 r_2 (\zeta^{6 \nu} z_2 z_4 z_5+\zeta^{2 \nu} z_3 z_5 z_6+\zeta^{5 \nu} z_1 z_4 z_6)+\\
  &-r_3 (\zeta^{10 \nu} z_1 z_2 z_6+\zeta^{4 \nu} z_1 z_3 z_5+\zeta^{12 \nu} z_2 z_3 z_4)+\\
  &+r_1 (\zeta^{3 \nu} z_3 z_4 z_5+\zeta^{9 \nu} z_2 z_4 z_6+\zeta^{\nu} z_1 z_5 z_6),
\endaligned$$
where
$$r_0=2(\theta_1-\theta_3)-3(\theta_2-\theta_4), \quad
  r_{\infty}=2(\theta_4-\theta_2)-3(\theta_1-\theta_3),$$
$$r_1=\sqrt{-13-2 \sqrt{13}}, \quad r_2=\sqrt{\frac{-13+3 \sqrt{13}}{2}},$$
$$r_3=\sqrt{-13+2 \sqrt{13}}, \quad r_4=\sqrt{\frac{-13-3 \sqrt{13}}{2}}.$$
This leads us to define the following senary cubic forms (cubic forms in six variables):
$$\left\{\aligned
  \mathbf{D}_0 &=z_1 z_2 z_3,\\
  \mathbf{D}_1 &=2 z_2 z_3^2+z_2^2 z_6-z_4^2 z_5+z_1 z_5 z_6,\\
  \mathbf{D}_2 &=-z_6^3+z_2^2 z_4-2 z_2 z_5^2+z_1 z_4 z_5+3 z_3 z_5 z_6,\\
  \mathbf{D}_3 &=2 z_1 z_2^2+z_1^2 z_5-z_4 z_6^2+z_3 z_4 z_5,\\
  \mathbf{D}_4 &=-z_2^2 z_3+z_1 z_6^2-2 z_4^2 z_6-z_1 z_3 z_5,\\
  \mathbf{D}_5 &=-z_4^3+z_3^2 z_5-2 z_3 z_6^2+z_2 z_5 z_6+3 z_1 z_4 z_6,\\
  \mathbf{D}_6 &=-z_5^3+z_1^2 z_6-2 z_1 z_4^2+z_3 z_4 z_6+3 z_2 z_4 z_5,\\
  \mathbf{D}_7 &=-z_2^3+z_3 z_4^2-z_1 z_3 z_6-3 z_1 z_2 z_5+2 z_1^2 z_4,\\
  \mathbf{D}_8 &=-z_1^3+z_2 z_6^2-z_2 z_3 z_5-3 z_1 z_3 z_4+2 z_3^2 z_6,\\
  \mathbf{D}_9 &=2 z_1^2 z_3+z_3^2 z_4-z_5^2 z_6+z_2 z_4 z_6,\\
  \mathbf{D}_{10} &=-z_1 z_3^2+z_2 z_4^2-2 z_4 z_5^2-z_1 z_2 z_6,\\
  \mathbf{D}_{11} &=-z_3^3+z_1 z_5^2-z_1 z_2 z_4-3 z_2 z_3 z_6+2 z_2^2 z_5,\\
  \mathbf{D}_{12} &=-z_1^2 z_2+z_3 z_5^2-2 z_5 z_6^2-z_2 z_3 z_4,\\
  \mathbf{D}_{\infty}&=z_4 z_5 z_6.
\endaligned\right.\eqno{(3.10)}$$
Then
$$\aligned
  &-13 \sqrt{13} ST^{\nu}(\mathbf{D}_0)\\
 =&r_0 \mathbf{D}_0+r_1 \zeta^{\nu} \mathbf{D}_1+
   r_2 \zeta^{2 \nu} \mathbf{D}_2+
   r_1 \zeta^{3 \nu} \mathbf{D}_3+r_3 \zeta^{4 \nu} \mathbf{D}_4+\\
  &+r_2 \zeta^{5 \nu} \mathbf{D}_5+r_2 \zeta^{6 \nu} \mathbf{D}_6+
   r_4 \zeta^{7 \nu} \mathbf{D}_7+r_4 \zeta^{8 \nu} \mathbf{D}_8+\\
  &+r_1 \zeta^{9 \nu} \mathbf{D}_9+r_3 \zeta^{10 \nu} \mathbf{D}_{10}
   +r_4 \zeta^{11 \nu} \mathbf{D}_{11}+r_3 \zeta^{12 \nu} \mathbf{D}_{12}
   +r_{\infty} \mathbf{D}_{\infty}.
\endaligned$$
$$\aligned
  &-13 \sqrt{13} ST^{\nu}(\mathbf{D}_{\infty})\\
 =&r_{\infty} \mathbf{D}_0-r_3 \zeta^{\nu} \mathbf{D}_1-
   r_4 \zeta^{2 \nu} \mathbf{D}_2-r_3 \zeta^{3 \nu} \mathbf{D}_3+
   r_1 \zeta^{4 \nu} \mathbf{D}_4+\\
  &-r_4 \zeta^{5 \nu} \mathbf{D}_5-r_4 \zeta^{6 \nu} \mathbf{D}_6+
   r_2 \zeta^{7 \nu} \mathbf{D}_7+r_2 \zeta^{8 \nu} \mathbf{D}_8+\\
  &-r_3 \zeta^{9 \nu} \mathbf{D}_9+r_1 \zeta^{10 \nu} \mathbf{D}_{10}+
   r_2 \zeta^{11 \nu} \mathbf{D}_{11}+r_1 \zeta^{12 \nu} \mathbf{D}_{12}-
   r_0 \mathbf{D}_{\infty}.
\endaligned$$

  Let
$$\delta_{\infty}(z_1, z_2, z_3, z_4, z_5, z_6)
 =13^2 (z_1^2 z_2^2 z_3^2+z_4^2 z_5^2 z_6^2)\eqno{(3.11)}$$
and
$$\delta_{\nu}(z_1, z_2, z_3, z_4, z_5, z_6)
 =\delta_{\infty}(ST^{\nu}(z_1, z_2, z_3, z_4, z_5, z_6))\eqno{(3.12)}$$
for $\nu=0, 1, \ldots, 12$. Then
$$\delta_{\nu}=13^2 ST^{\nu}(\mathbf{G}_0)
 =-13 \mathbf{G}_0+\zeta^{\nu} \mathbf{G}_1+\zeta^{2 \nu} \mathbf{G}_2+
  \cdots+\zeta^{12 \nu} \mathbf{G}_{12},\eqno{(3.13)}$$
where the senary sextic forms (i.e., sextic forms in six
variables) are given as follows:
$$\left\{\aligned
  \mathbf{G}_0 =&\mathbf{D}_0^2+\mathbf{D}_{\infty}^2,\\
  \mathbf{G}_1 =&-\mathbf{D}_7^2+2 \mathbf{D}_0 \mathbf{D}_1+10 \mathbf{D}_{\infty}
                 \mathbf{D}_1+2 \mathbf{D}_2 \mathbf{D}_{12}+\\
                &-2 \mathbf{D}_3 \mathbf{D}_{11}-4 \mathbf{D}_4 \mathbf{D}_{10}-2
                 \mathbf{D}_9 \mathbf{D}_5,\\
  \mathbf{G}_2 =&-2 \mathbf{D}_1^2-4 \mathbf{D}_0 \mathbf{D}_2+6 \mathbf{D}_{\infty}
                 \mathbf{D}_2-2 \mathbf{D}_4 \mathbf{D}_{11}+\\
                &+2 \mathbf{D}_5 \mathbf{D}_{10}-2 \mathbf{D}_6 \mathbf{D}_9-2
                 \mathbf{D}_7 \mathbf{D}_8,\\
  \mathbf{G}_3 =&-\mathbf{D}_8^2+2 \mathbf{D}_0 \mathbf{D}_3+10 \mathbf{D}_{\infty}
                 \mathbf{D}_3+2 \mathbf{D}_6 \mathbf{D}_{10}+\\
                &-2 \mathbf{D}_9 \mathbf{D}_7-4 \mathbf{D}_{12} \mathbf{D}_4-2
                 \mathbf{D}_1 \mathbf{D}_2,\\
  \mathbf{G}_4 =&-\mathbf{D}_2^2+10 \mathbf{D}_0 \mathbf{D}_4-2 \mathbf{D}_{\infty}
                 \mathbf{D}_4+2 \mathbf{D}_5 \mathbf{D}_{12}+\\
                &-2 \mathbf{D}_9 \mathbf{D}_8-4 \mathbf{D}_1 \mathbf{D}_3-2
                 \mathbf{D}_{10} \mathbf{D}_7,\\
  \mathbf{G}_5 =&-2 \mathbf{D}_9^2-4 \mathbf{D}_0 \mathbf{D}_5+6 \mathbf{D}_{\infty}
                 \mathbf{D}_5-2 \mathbf{D}_{10} \mathbf{D}_8+\\
                &+2 \mathbf{D}_6 \mathbf{D}_{12}-2 \mathbf{D}_2 \mathbf{D}_3-2
                 \mathbf{D}_{11} \mathbf{D}_7,\\
  \mathbf{G}_6 =&-2 \mathbf{D}_3^2-4 \mathbf{D}_0 \mathbf{D}_6+6 \mathbf{D}_{\infty}
                 \mathbf{D}_6-2 \mathbf{D}_{12} \mathbf{}_7+\\
                &+2 \mathbf{D}_2 \mathbf{D}_4-2 \mathbf{D}_5 \mathbf{D}_1-2
                 \mathbf{D}_8 \mathbf{D}_{11},\\
  \mathbf{G}_7 =&-2 \mathbf{D}_{10}^2+6 \mathbf{D}_0 \mathbf{D}_7+4 \mathbf{D}_{\infty}
                 \mathbf{D}_7-2 \mathbf{D}_1 \mathbf{D}_6+\\
                &-2 \mathbf{D}_2 \mathbf{D}_5-2 \mathbf{D}_8 \mathbf{D}_{12}-2
                 \mathbf{D}_9 \mathbf{D}_{11},\\
  \mathbf{G}_8 =&-2 \mathbf{D}_4^2+6 \mathbf{D}_0 \mathbf{D}_8+4 \mathbf{D}_{\infty}
                 \mathbf{D}_8-2 \mathbf{D}_3 \mathbf{D}_5+\\
                &-2 \mathbf{D}_6 \mathbf{D}_2-2 \mathbf{D}_{11} \mathbf{D}_{10}-2
                 \mathbf{D}_1 \mathbf{D}_7,\\
  \mathbf{G}_9 =&-\mathbf{D}_{11}^2+2 \mathbf{D}_0 \mathbf{D}_9+10 \mathbf{D}_{\infty}
                 \mathbf{D}_9+2 \mathbf{D}_5 \mathbf{D}_4+\\
                &-2 \mathbf{D}_1 \mathbf{D}_8-4 \mathbf{D}_{10} \mathbf{D}_{12}-2
                 \mathbf{D}_3 \mathbf{D}_6,\\
  \mathbf{G}_{10} =&-\mathbf{D}_5^2+10 \mathbf{D}_0 \mathbf{D}_{10}-2 \mathbf{D}_{\infty}
                    \mathbf{D}_{10}+2 \mathbf{D}_6 \mathbf{D}_4+\\
                   &-2 \mathbf{D}_3 \mathbf{D}_7-4 \mathbf{D}_9 \mathbf{D}_1-2
                    \mathbf{D}_{12} \mathbf{D}_{11},\\
  \mathbf{G}_{11} =&-2 \mathbf{D}_{12}^2+6 \mathbf{D}_0 \mathbf{D}_{11}+4 \mathbf{D}_{\infty}
                    \mathbf{D}_{11}-2 \mathbf{D}_9 \mathbf{D}_2+\\
                   &-2 \mathbf{D}_5 \mathbf{D}_6-2 \mathbf{D}_7 \mathbf{D}_4-2
                    \mathbf{D}_3 \mathbf{D}_8,\\
  \mathbf{G}_{12} =&-\mathbf{D}_6^2+10 \mathbf{D}_0 \mathbf{D}_{12}-2 \mathbf{D}_{\infty}
                    \mathbf{D}_{12}+2 \mathbf{D}_2 \mathbf{D}_{10}+\\
                   &-2 \mathbf{D}_1 \mathbf{D}_{11}-4 \mathbf{D}_3 \mathbf{D}_9-2
                    \mathbf{D}_4 \mathbf{D}_8.
\endaligned\right.\eqno{(3.14)}$$
We have that $\mathbf{G}_0$ is invariant under the action of
$\langle H, T \rangle$, a maximal subgroup of order $78$ of $G$
with index $14$. Note that $\delta_{\infty}$, $\delta_{\nu}$ for
$\nu=0, \ldots, 12$ form an algebraic equation of degree fourteen.
However, we have $\delta_{\infty}+\sum_{\nu=0}^{12}
\delta_{\nu}=0$. Hence, it is not the Jacobian equation of degree
fourteen.

  Recall that the theta functions with characteristic
$\begin{bmatrix} \epsilon\\ \epsilon^{\prime} \end{bmatrix} \in
\mathbb{R}^2$ is defined by the following series which converges
uniformly and absolutely on compact subsets of $\mathbb{C} \times
\mathbb{H}$ (see \cite{FK}, p. 73):
$$\theta \begin{bmatrix} \epsilon\\ \epsilon^{\prime} \end{bmatrix}(z, \tau)
 =\sum_{n \in \mathbb{Z}} \exp \left\{2 \pi i \left[\frac{1}{2}
  \left(n+\frac{\epsilon}{2}\right)^2 \tau+\left(n+\frac{\epsilon}{2}\right)
  \left(z+\frac{\epsilon^{\prime}}{2}\right)\right]\right\}.$$
The modified theta constants (see \cite{FK}, p. 215)
$\varphi_l(\tau):=\theta [\chi_l](0, k \tau)$,
where the characteristic $\chi_l=\begin{bmatrix} \frac{2l+1}{k}\\ 1
\end{bmatrix}$, $l=0, \ldots, \frac{k-3}{2}$, for odd $k$ and
$\chi_l=\begin{bmatrix} \frac{2l}{k}\\ 0 \end{bmatrix}$, $l=0,
\ldots, \frac{k}{2}$, for even $k$. We have the following:

\textbf{Proposition 3.1.} (see \cite{FK}, p. 236). {\it For each odd
integer $k \geq 5$, the map
$\Phi: \tau \mapsto (\varphi_0(\tau), \varphi_1(\tau), \ldots,
 \varphi_{\frac{k-5}{2}}(\tau), \varphi_{\frac{k-3}{2}}(\tau))$
from $\mathbb{H} \cup \mathbb{Q} \cup \{ \infty \}$ to
$\mathbb{C}^{\frac{k-1}{2}}$, defines a holomorphic mapping from
$\overline{\mathbb{H}/\Gamma(k)}$ into $\mathbb{C} \mathbb{P}^{\frac{k-3}{2}}$.}

  In our case, the map
$\Phi: \tau \mapsto (\varphi_0(\tau), \varphi_1(\tau), \varphi_2(\tau),
 \varphi_3(\tau), \varphi_4(\tau), \varphi_5(\tau))$
gives a holomorphic mapping from the modular curve
$X(13)=\overline{\mathbb{H}/\Gamma(13)}$ into $\mathbb{C} \mathbb{P}^5$,
which corresponds to our six-dimensional representation, i.e., up
to the constants, $z_1, \ldots, z_6$ are just modular forms
$\varphi_0(\tau), \ldots, \varphi_5(\tau)$. Let
$$\left\{\aligned
  a_1(z) &:=e^{-\frac{11 \pi i}{26}} \theta
            \begin{bmatrix} \frac{11}{13}\\ 1 \end{bmatrix}(0, 13z)
           =q^{\frac{121}{104}} \sum_{n \in \mathbb{Z}} (-1)^n q^{\frac{1}{2}(13n^2+11n)},\\
  a_2(z) &:=e^{-\frac{7 \pi i}{26}} \theta
            \begin{bmatrix} \frac{7}{13}\\ 1 \end{bmatrix}(0, 13z)
           =q^{\frac{49}{104}} \sum_{n \in \mathbb{Z}} (-1)^n q^{\frac{1}{2}(13n^2+7n)},\\
  a_3(z) &:=e^{-\frac{5 \pi i}{26}} \theta
            \begin{bmatrix} \frac{5}{13}\\ 1 \end{bmatrix}(0, 13z)
           =q^{\frac{25}{104}} \sum_{n \in \mathbb{Z}} (-1)^n q^{\frac{1}{2}(13n^2+5n)},\\
  a_4(z) &:=-e^{-\frac{3 \pi i}{26}} \theta
            \begin{bmatrix} \frac{3}{13}\\ 1 \end{bmatrix}(0, 13z)
           =-q^{\frac{9}{104}} \sum_{n \in \mathbb{Z}} (-1)^n q^{\frac{1}{2}(13n^2+3n)},\\
  a_5(z) &:=e^{-\frac{9 \pi i}{26}} \theta
            \begin{bmatrix} \frac{9}{13}\\ 1 \end{bmatrix}(0, 13z)
           =q^{\frac{81}{104}} \sum_{n \in \mathbb{Z}} (-1)^n q^{\frac{1}{2}(13n^2+9n)},\\
  a_6(z) &:=e^{-\frac{\pi i}{26}} \theta
            \begin{bmatrix} \frac{1}{13}\\ 1 \end{bmatrix}(0, 13z)
           =q^{\frac{1}{104}} \sum_{n \in \mathbb{Z}} (-1)^n q^{\frac{1}{2}(13n^2+n)}
\endaligned\right.\eqno{(3.15)}$$
be the theta constants of order $13$ and
$$\mathbf{A}(z):=(a_1(z), a_2(z), a_3(z), a_4(z), a_5(z), a_6(z))^{T}.$$
The significance of our six dimensional representation of $\text{PSL}(2, 13)$
comes from the following:

\textbf{Proposition 3.2} (see \cite{Y2}, Proposition 2.5). {\it If $z \in \mathbb{H}$,
then the following relations hold:
$$\mathbf{A}(z+1)=e^{-\frac{3 \pi i}{4}} T \mathbf{A}(z), \quad
  \mathbf{A}\left(-\frac{1}{z}\right)=e^{\frac{\pi i}{4}} \sqrt{z}
  S \mathbf{A}(z),\eqno{(3.16)}$$
where $S$ and $T$ are given in (3.1) and (3.2), and
$0<\text{arg} \sqrt{z} \leq \pi/2$.}

  Recall that the principal congruence subgroup of level $13$ is the
normal subgroup $\Gamma(13)$ of $\Gamma=\text{PSL}(2, \mathbb{Z})$
defined by the exact sequence
$1 \rightarrow \Gamma(13) \rightarrow \Gamma(1) \stackrel{f}{\rightarrow} G \rightarrow 1$
where $f(\gamma) \equiv \gamma$ (mod $13$) for $\gamma \in \Gamma=\Gamma(1)$.
There is a representation $\rho: \Gamma \rightarrow \text{PGL}(6, \mathbb{C})$
with kernel $\Gamma(13)$ defined as follows: if
$t=\begin{pmatrix} 1 & 1\\ 0 & 1 \end{pmatrix}$ and
$s=\begin{pmatrix} 0 & -1\\ 1 & 0 \end{pmatrix}$,
then $\rho(t)=T$ and $\rho(s)=S$. To see that such a representation
exists, note that $\Gamma$ is defined by the presentation
$\langle s, t; s^2=(st)^3=1 \rangle$ satisfied by $s$ and $t$ and we
have proved that $S$ and $T$ satisfy these relations. Moreover, we
have proved that $G$ is defined by the presentation
$\langle S, T; S^2=T^{13}=(ST)^3=1 \rangle$. Let $p=s t^{-1} s$
and $q=s t^3$. Then
$$h:=q^5 p^2 \cdot p^2 q^6 p^8 \cdot q^5 p^2 \cdot p^3 q
    =\begin{pmatrix}
     4, 428, 249 & -10, 547, 030\\
    -11, 594, 791 & 27, 616, 019
     \end{pmatrix}$$
satisfies that $\rho(h)=H$. The off-diagonal elements of the
matrix $h$, which corresponds to $H$, are congruent to $0$
mod $13$. The connection to $\Gamma_0(13)$ should be obvious.

  Put $x_i(z)=\eta(z) a_i(z)$ and $y_i(z)=\eta^3(z) a_i(z)$
$(1 \leq i \leq 6)$. Let
$$X(z)=(x_1(z), \ldots, x_6(z))^{T} \quad \text{and} \quad
  Y(z)=(y_1(z), \ldots, y_6(z))^{T}.$$
Then
$X(z)=\eta(z) \mathbf{A}(z)$ and $Y(z)=\eta^3(z) \mathbf{A}(z)$.
Recall that $\eta(z)$ satisfies the following transformation formulas
$\eta(z+1)=e^{\frac{\pi i}{12}} \eta(z)$ and
$\eta\left(-\frac{1}{z}\right)=e^{-\frac{\pi i}{4}} \sqrt{z} \eta(z)$.
By Proposition 3.2, we have
$$X(z+1)=e^{-\frac{2 \pi i}{3}} \rho(t) X(z), \quad
  X\left(-\frac{1}{z}\right)=z \rho(s) X(z),$$
$$Y(z+1)=e^{-\frac{\pi i}{2}} \rho(t) Y(z), \quad
  Y\left(-\frac{1}{z}\right)=e^{-\frac{\pi i}{2}} z^2 \rho(s) Y(z).$$
Define $j(\gamma, z):=cz+d$ if $z \in \mathbb{H}$ and
$\gamma=\begin{pmatrix} a & b\\ c & d \end{pmatrix} \in \Gamma(1)$.
Hence, $X(\gamma(z))=u(\gamma) j(\gamma, z) \rho(\gamma) X(z)$ and
$Y(\gamma(z))=v(\gamma) j(\gamma, z)^2 \rho(\gamma) Y(z)$ for
$\gamma \in \Gamma(1)$, where $u(\gamma)=1, \omega$ or $\omega^2$
with $\omega=e^{\frac{2 \pi i}{3}}$ and $v(\gamma)=\pm 1$ or $\pm i$.
Since $\Gamma(13)=\text{ker}$ $\rho$, we have
$X(\gamma(z))=u(\gamma) j(\gamma, z) X(z)$ and
$Y(\gamma(z))=v(\gamma) j(\gamma, z)^2 Y(z)$ for $\gamma \in \Gamma(13)$.
This means that the functions $x_1(z)$, $\ldots$, $x_6(z)$ are modular
forms of weight one for $\Gamma(13)$ with the same multiplier
$u(\gamma)=1, \omega$ or $\omega^2$ and $y_1(z)$, $\ldots$, $y_6(z)$
are modular forms of weight two for $\Gamma(13)$ with the
same multiplier $v(\gamma)=\pm 1 $ or $\pm i$.

  From now on, we will use the following abbreviation
$$\mathbf{A}_j=\mathbf{A}_j(a_1(z), \ldots, a_6(z)) \quad
  (0 \leq j \leq 6),$$
$$\mathbf{D}_j=\mathbf{D}_j(a_1(z), \ldots, a_6(z)) \quad
  (j=0,1, \ldots, 12, \infty)$$
and
$$\mathbf{G}_j=\mathbf{G}_j(a_1(z), \ldots, a_6(z)) \quad
  (0 \leq j \leq 12).$$
We have
$$\left\{\aligned
  \mathbf{A}_0 &=q^{\frac{1}{4}} (1+O(q)),\\
  \mathbf{A}_1 &=q^{\frac{17}{52}} (2+O(q)),\\
  \mathbf{A}_2 &=q^{\frac{29}{52}} (2+O(q)),\\
  \mathbf{A}_3 &=q^{\frac{49}{52}} (1+O(q)),\\
  \mathbf{A}_4 &=q^{\frac{25}{52}} (-1+O(q)),\\
  \mathbf{A}_5 &=q^{\frac{9}{52}} (-1+O(q)),\\
  \mathbf{A}_6 &=q^{\frac{1}{52}} (-1+O(q)),
\endaligned\right.$$
and
$$\left\{\aligned
  \mathbf{D}_0 &=q^{\frac{15}{8}} (1+O(q)),\\
  \mathbf{D}_{\infty} &=q^{\frac{7}{8}} (-1+O(q)),\\
  \mathbf{D}_1 &=q^{\frac{99}{104}} (2+O(q)),\\
  \mathbf{D}_2 &=q^{\frac{3}{104}} (-1+O(q)),\\
  \mathbf{D}_3 &=q^{\frac{11}{104}} (1+O(q)),\\
  \mathbf{D}_4 &=q^{\frac{19}{104}} (-2+O(q)),\\
  \mathbf{D}_5 &=q^{\frac{27}{104}} (-1+O(q)),
\endaligned\right. \quad \quad
  \left\{\aligned
  \mathbf{D}_6 &=q^{\frac{35}{104}} (-1+O(q)),\\
  \mathbf{D}_7 &=q^{\frac{43}{104}} (1+O(q)),\\
  \mathbf{D}_8 &=q^{\frac{51}{104}} (3+O(q)),\\
  \mathbf{D}_9 &=q^{\frac{59}{104}} (-2+O(q)),\\
  \mathbf{D}_{10} &=q^{\frac{67}{104}} (1+O(q)),\\
  \mathbf{D}_{11} &=q^{\frac{75}{104}} (-4+O(q)),\\
  \mathbf{D}_{12} &=q^{\frac{83}{104}} (-1+O(q)).
\endaligned\right.$$
Hence,
$$\left\{\aligned
  \mathbf{G}_0 &=q^{\frac{7}{4}} (1+O(q)),\\
  \mathbf{G}_1 &=q^{\frac{43}{52}} (13+O(q)),\\
  \mathbf{G}_2 &=q^{\frac{47}{52}} (-22+O(q)),\\
  \mathbf{G}_3 &=q^{\frac{51}{52}} (-21+O(q)),\\
  \mathbf{G}_4 &=q^{\frac{3}{52}} (-1+O(q)),\\
  \mathbf{G}_5 &=q^{\frac{7}{52}} (2+O(q)),\\
  \mathbf{G}_6 &=q^{\frac{11}{52}} (2+O(q)),
\endaligned\right. \quad \quad
  \left\{\aligned
  \mathbf{G}_7 &=q^{\frac{15}{52}} (-2+O(q)),\\
  \mathbf{G}_8 &=q^{\frac{19}{52}} (-8+O(q)),\\
  \mathbf{G}_9 &=q^{\frac{23}{52}} (6+O(q)),\\
  \mathbf{G}_{10} &=q^{\frac{27}{52}} (1+O(q)),\\
  \mathbf{G}_{11} &=q^{\frac{31}{52}} (-8+O(q)),\\
  \mathbf{G}_{12} &=q^{\frac{35}{52}} (17+O(q)).
\endaligned\right.$$
Note that
$$\aligned
  w_{\nu} &=(\mathbf{A}_0+\zeta^{\nu} \mathbf{A}_1+\zeta^{4 \nu} \mathbf{A}_2
           +\zeta^{9 \nu} \mathbf{A}_3+\zeta^{3 \nu} \mathbf{A}_4
           +\zeta^{12 \nu} \mathbf{A}_5+\zeta^{10 \nu} \mathbf{A}_6)^2\\
          &=\mathbf{A}_0^2+2 (\mathbf{A}_1 \mathbf{A}_5+\mathbf{A}_2 \mathbf{A}_3
           +\mathbf{A}_4 \mathbf{A}_6)+\\
          &+2 \zeta^{\nu} (\mathbf{A}_0 \mathbf{A}_1+\mathbf{A}_2 \mathbf{A}_6)
           +2 \zeta^{3 \nu} (\mathbf{A}_0 \mathbf{A}_4+\mathbf{A}_2 \mathbf{A}_5)+\\
          &+2 \zeta^{9 \nu} (\mathbf{A}_0 \mathbf{A}_3+\mathbf{A}_5 \mathbf{A}_6)
           +2 \zeta^{12 \nu} (\mathbf{A}_0 \mathbf{A}_5+\mathbf{A}_3 \mathbf{A}_4)+\\
          &+2 \zeta^{10 \nu} (\mathbf{A}_0 \mathbf{A}_6+\mathbf{A}_1 \mathbf{A}_3)
           +2 \zeta^{4 \nu} (\mathbf{A}_0 \mathbf{A}_2+\mathbf{A}_1 \mathbf{A}_4)+\\
          &+\zeta^{2 \nu} (\mathbf{A}_1^2+2 \mathbf{A}_4 \mathbf{A}_5)
           +\zeta^{5 \nu} (\mathbf{A}_3^2+2 \mathbf{A}_1 \mathbf{A}_2)+\\
          &+\zeta^{6 \nu} (\mathbf{A}_4^2+2 \mathbf{A}_3 \mathbf{A}_6)
           +\zeta^{11 \nu} (\mathbf{A}_5^2+2 \mathbf{A}_1 \mathbf{A}_6)+\\
          &+\zeta^{8 \nu} (\mathbf{A}_2^2+2 \mathbf{A}_3 \mathbf{A}_5)
           +\zeta^{7 \nu} (\mathbf{A}_6^2+2 \mathbf{A}_4 \mathbf{A}_2),
\endaligned$$
where
$$\mathbf{A}_0^2+2 (\mathbf{A}_1 \mathbf{A}_5+\mathbf{A}_2 \mathbf{A}_3
  +\mathbf{A}_4 \mathbf{A}_6)=q^{\frac{1}{2}} (-1+O(q)),$$
$$\left\{\aligned
  \mathbf{A}_0 \mathbf{A}_1+\mathbf{A}_2 \mathbf{A}_6 &=q^{\frac{41}{26}} (-3+O(q)),\\
  \mathbf{A}_0 \mathbf{A}_4+\mathbf{A}_2 \mathbf{A}_5 &=q^{\frac{19}{26}} (-3+O(q)),\\
  \mathbf{A}_0 \mathbf{A}_3+\mathbf{A}_5 \mathbf{A}_6 &=q^{\frac{5}{26}} (1+O(q)),\\
  \mathbf{A}_0 \mathbf{A}_5+\mathbf{A}_3 \mathbf{A}_4 &=q^{\frac{11}{26}} (-1+O(q)),\\
  \mathbf{A}_0 \mathbf{A}_6+\mathbf{A}_1 \mathbf{A}_3 &=q^{\frac{7}{26}} (-1+O(q)),\\
  \mathbf{A}_0 \mathbf{A}_2+\mathbf{A}_1 \mathbf{A}_4 &=q^{\frac{47}{26}} (-1+O(q)),
\endaligned\right.$$
and
$$\left\{\aligned
  \mathbf{A}_1^2+2 \mathbf{A}_4 \mathbf{A}_5 &=q^{\frac{17}{26}} (6+O(q)),\\
  \mathbf{A}_3^2+2 \mathbf{A}_1 \mathbf{A}_2 &=q^{\frac{23}{26}} (8+O(q)),\\
  \mathbf{A}_4^2+2 \mathbf{A}_3 \mathbf{A}_6 &=q^{\frac{25}{26}} (-1+O(q)),\\
  \mathbf{A}_5^2+2 \mathbf{A}_1 \mathbf{A}_6 &=q^{\frac{9}{26}} (-3+O(q)),\\
  \mathbf{A}_2^2+2 \mathbf{A}_3 \mathbf{A}_5 &=q^{\frac{29}{26}} (2+O(q)),\\
  \mathbf{A}_6^2+2 \mathbf{A}_4 \mathbf{A}_2 &=q^{\frac{1}{26}} (1+O(q)).
\endaligned\right.$$

{\it Proof of Theorem 1.1}. We divide the proof into three parts (see also
\cite{Y3}). The first part is the calculation of $\Phi_{20}$ and
$\Phi_{12}^{\prime}$. Let
$$\Phi_{20}=w_0^5+w_1^5+\cdots+w_{12}^5+w_{\infty}^5.$$
As a polynomial in six variables, $\Phi_{20}(z_1, z_2, z_3, z_4, z_5, z_6)$
is a $G$-invariant polynomial. Moreover, for $\gamma \in \Gamma(1)$,
$$\aligned
  &\Phi_{20}(Y(\gamma(z))^{T})=\Phi_{20}(v(\gamma)
   j(\gamma, z)^2 (\rho(\gamma) Y(z))^{T})\\
 =&v(\gamma)^{20} j(\gamma, z)^{40} \Phi_{20}((\rho(\gamma) Y(z))^{T})
 =j(\gamma, z)^{40} \Phi_{20}((\rho(\gamma) Y(z))^{T}).
\endaligned$$
Note that $\rho(\gamma) \in \langle \rho(s), \rho(t) \rangle=G$ and
$\Phi_{20}$ is a $G$-invariant polynomial, we have
$$\Phi_{20}(Y(\gamma(z))^{T})=j(\gamma, z)^{40} \Phi_{20}(Y(z)^{T}),
  \quad \text{for $\gamma \in \Gamma(1)$}.$$
This implies that $\Phi_{20}(y_1(z), \ldots, y_6(z))$ is a modular form
of weight $40$ for the full modular group $\Gamma(1)$. Moreover, we will
show that it is a cusp form. In fact,
$$\aligned
  &\Phi_{20}(a_1(z), \ldots, a_6(z))=13^5 q^{\frac{5}{2}} (1+O(q))^5+\\
  &+\sum_{\nu=0}^{12} [q^{\frac{1}{2}} (-1+O(q))+\\
  &+2 \zeta^{\nu} q^{\frac{41}{26}} (-3+O(q))+2 \zeta^{3 \nu} q^{\frac{19}{26}} (-3+O(q))
   +2 \zeta^{9 \nu} q^{\frac{5}{26}} (1+O(q))+\\
  &+2 \zeta^{12 \nu} q^{\frac{11}{26}} (-1+O(q))+2 \zeta^{10 \nu} q^{\frac{7}{26}} (-1+O(q))
   +2 \zeta^{4 \nu} q^{\frac{47}{26}} (-1+O(q))+\\
  &+\zeta^{2 \nu} q^{\frac{17}{26}} (6+O(q))+\zeta^{5 \nu} q^{\frac{23}{26}} (8+O(q))
   +\zeta^{6 \nu} q^{\frac{25}{26}} (-1+O(q))+\\
  &+\zeta^{11 \nu} q^{\frac{9}{26}} (-3+O(q))+\zeta^{8 \nu} q^{\frac{29}{26}} (2+O(q))
   +\zeta^{7 \nu} q^{\frac{1}{26}} (1+O(q))]^5.
\endaligned$$
We will calculate the $q^{\frac{1}{2}}$-term which is the lowest degree.
For the partition $13=4 \cdot 1+9$, the corresponding term is
$$\begin{pmatrix} 5\\ 4, 1 \end{pmatrix} (\zeta^{7 \nu} q^{\frac{1}{26}})^4
  \cdot (-3) \zeta^{11 \nu} q^{\frac{9}{26}}=-15 q^{\frac{1}{2}}.$$
For the partition $13=3 \cdot 1+2 \cdot 5$, the corresponding term is
$$\begin{pmatrix} 5\\ 3, 2 \end{pmatrix} (\zeta^{7 \nu} q^{\frac{1}{26}})^3
  \cdot (2 \zeta^{9 \nu} q^{\frac{5}{26}})^2=40 q^{\frac{1}{2}}.$$
Hence, for $\Phi_{20}(y_1(z), \ldots, y_6(z))$ which is a modular form for
$\Gamma(1)$ with weight $40$, the lowest degree term is given by
$$(-15+40) q^{\frac{1}{2}} \cdot q^{\frac{3}{24} \cdot 20}=25 q^3.$$
Thus,
$$\Phi_{20}(y_1(z), \ldots, y_6(z))=q^3 (13 \cdot 25+O(q)).$$
The leading term of $\Phi_{20}(y_1(z), \ldots, y_6(z))$ together with
its weight $40$ suffice to identify this modular form with
$\Phi_{20}(y_1(z), \ldots, y_6(z))=13 \cdot 25 \Delta(z)^3 E_4(z)$.
Consequently,
$$\Phi_{20}(x_1(z), \ldots, x_6(z))=13 \cdot 25 \Delta(z)^3
  E_4(z)/\eta(z)^{40}=13 \cdot 25 \eta(z)^8 \Delta(z) E_4(z).$$

  Let
$$\Phi_{12}^{\prime}=w_0^3+w_1^3+\cdots+w_{12}^3+w_{\infty}^3.$$
The calculation of $\Phi_{12}^{\prime}$ is similar as that of $\Phi_{20}$.
We find that
$$\Phi_{12}^{\prime}(x_1(z), \ldots, x_6(z))=-13 \cdot 30 \Delta(z).$$

  The second part is the calculation of $\Phi_4$, $\Phi_8$ and $\Phi_{16}$.
The calculation of $\Phi_4$ has been done in \cite{Y2}, Theorem 3.1. We
will give the calculation of $\Phi_{16}$. Let
$$\Phi_{16}=w_0^4+w_1^4+\cdots+w_{12}^4+w_{\infty}^4.$$
Similar as the above calculation for $\Phi_{20}$, we find that
$\Phi_{16}(y_1(z), \ldots, y_6(z))$ is a modular form of weight $32$ for
the full modular group $\Gamma(1)$. Moreover, we will show that it is a
cusp form. In fact,
$$\aligned
  &\Phi_{16}(a_1(z), \ldots, a_6(z))=13^4 q^2 (1+O(q))^4+\\
  &+\sum_{\nu=0}^{12} [q^{\frac{1}{2}} (-1+O(q))+\\
  &+2 \zeta^{\nu} q^{\frac{41}{26}} (-3+O(q))+2 \zeta^{3 \nu} q^{\frac{19}{26}} (-3+O(q))
   +2 \zeta^{9 \nu} q^{\frac{5}{26}} (1+O(q))+\\
  &+2 \zeta^{12 \nu} q^{\frac{11}{26}} (-1+O(q))+2 \zeta^{10 \nu} q^{\frac{7}{26}} (-1+O(q))
   +2 \zeta^{4 \nu} q^{\frac{47}{26}} (-1+O(q))+\\
  &+\zeta^{2 \nu} q^{\frac{17}{26}} (6+O(q))+\zeta^{5 \nu} q^{\frac{23}{26}} (8+O(q))
   +\zeta^{6 \nu} q^{\frac{25}{26}} (-1+O(q))+\\
  &+\zeta^{11 \nu} q^{\frac{9}{26}} (-3+O(q))+\zeta^{8 \nu} q^{\frac{29}{26}} (2+O(q))
   +\zeta^{7 \nu} q^{\frac{1}{26}} (1+O(q))]^4.
\endaligned$$
We will calculate the $q$-term which is the lowest degree. For example,
consider the partition $26=3 \cdot 1+23$, the corresponding term is
$$\begin{pmatrix} 4\\ 3, 1 \end{pmatrix} (\zeta^{7 \nu} q^{\frac{1}{26}})^3
  \cdot 8 \zeta^{5 \nu} q^{\frac{23}{26}}=32 q.$$
For the other partitions, the calculation is similar. In conclusion,
we find that the coefficients of the $q$-term is an integer. Hence,
for $\Phi_{16}(y_1(z), \ldots, y_6(z))$ which is a modular form for
$\Gamma(1)$ with weight $32$, the lowest degree term is given by
$$\text{some integer} \cdot q \cdot q^{\frac{3}{24} \cdot 16}
 =\text{some integer} \cdot q^3.$$
This implies that $\Phi_{16}(y_1(z), \ldots, y_6(z))$ has a factor of
$\Delta(z)^3$, which is a cusp form of weight $36$. Therefore,
$\Phi_{16}(y_1(z), \ldots, y_6(z))=0$. The calculation of $\Phi_8$ is
similar as that of $\Phi_{16}$.

  The third part is the calculation of $\Phi_{12}$, $\Phi_{18}$ and
$\Phi_{30}$. Let
$$\Phi_{12}=\delta_0^2+\delta_1^2+\cdots+\delta_{12}^2+\delta_{\infty}^2.$$
As a polynomial in six variables, $\Phi_{12}(z_1, z_2, z_3, z_4, z_5, z_6)$
is a $G$-invariant polynomial. Moreover, for $\gamma \in \Gamma(1)$,
$$\aligned
  &\Phi_{12}(X(\gamma(z))^{T})
 =\Phi_{12}(u(\gamma) j(\gamma, z) (\rho(\gamma) X(z))^{T})\\
 =&u(\gamma)^{12} j(\gamma, z)^{12} \Phi_{12}((\rho(\gamma) X(z))^{T})
 =j(\gamma, z)^{12} \Phi_{12}((\rho(\gamma) X(z))^{T}).
\endaligned$$
Note that $\rho(\gamma) \in \langle \rho(s), \rho(t) \rangle=G$
and $\Phi_{12}$ is a $G$-invariant polynomial, we have
$$\Phi_{12}(X(\gamma(z))^{T})=j(\gamma, z)^{12} \Phi_{12}(X(z)^{T}),
  \quad \text{for $\gamma \in \Gamma(1)$}.$$
This implies that $\Phi_{12}(x_1(z), \ldots, x_6(z))$ is a modular form
of weight $12$ for the full modular group $\Gamma(1)$. Moreover, we will
show that it is a cusp form. In fact,
$$\aligned
  &\Phi_{12}(a_1(z), \ldots, a_6(z))=13^4 q^{\frac{7}{2}} (1+O(q))^2+\\
  &+\sum_{\nu=0}^{12} [-13 q^{\frac{7}{4}} (1+O(q))+\\
  &+\zeta^{\nu} q^{\frac{43}{52}} (13+O(q))+\zeta^{2 \nu} q^{\frac{47}{52}} (-22+O(q))
   +\zeta^{3 \nu} q^{\frac{51}{52}} (-21+O(q))+\\
  &+\zeta^{4 \nu} q^{\frac{3}{52}} (-1+O(q))+\zeta^{5 \nu} q^{\frac{7}{52}} (2+O(q))
   +\zeta^{6 \nu} q^{\frac{11}{52}} (2+O(q))+\\
  &+\zeta^{7 \nu} q^{\frac{15}{52}} (-2+O(q))+\zeta^{8 \nu} q^{\frac{19}{52}} (-8+O(q))
   +\zeta^{9 \nu} q^{\frac{23}{52}} (6+O(q))+\\
  &+\zeta^{10 \nu} q^{\frac{27}{52}} (1+O(q))+\zeta^{11 \nu} q^{\frac{31}{52}} (-8+O(q))
   +\zeta^{12 \nu} q^{\frac{35}{52}} (17+O(q))]^2.
\endaligned$$
We will calculate the $q^{\frac{1}{2}}$-term which is the lowest degree.
For the partition $26=3+23$, the corresponding term is
$$\begin{pmatrix} 2\\ 1, 1 \end{pmatrix} \zeta^{4 \nu} q^{\frac{3}{52}}
  \cdot (-1) \cdot \zeta^{9 \nu} q^{\frac{23}{52}} \cdot 6=-12 q^{\frac{1}{2}}.$$
For the partition $26=7+19$, the corresponding term is
$$\begin{pmatrix} 2\\ 1, 1 \end{pmatrix} \zeta^{5 \nu} q^{\frac{7}{52}}
  \cdot 2 \cdot \zeta^{8 \nu} q^{\frac{19}{52}} \cdot (-8)=-32 q^{\frac{1}{2}}.$$
For the partition $26=11+15$, the corresponding term is
$$\begin{pmatrix} 2\\ 1, 1 \end{pmatrix} \zeta^{6 \nu} q^{\frac{11}{52}}
  \cdot 2 \cdot \zeta^{7 \nu} q^{\frac{15}{52}} \cdot (-2)=-8 q^{\frac{1}{2}}.$$
Hence, for $\Phi_{12}(x_1(z), \ldots, x_6(z))$ which is a modular form for
$\Gamma(1)$ with weight $12$, the lowest degree term is given by
$(-12-32-8) q^{\frac{1}{2}} \cdot q^{\frac{12}{24}}=-52 q$.
Thus,
$$\Phi_{12}(x_1(z), \ldots, x_6(z))=q (-13 \cdot 52+O(q)).$$
The leading term of $\Phi_{12}(x_1(z), \ldots, x_6(z))$ together with its
weight $12$ suffice to identify this modular form with
$$\Phi_{12}(x_1(z), \ldots, x_6(z))=-13 \cdot 52 \Delta(z).$$

  Let
$$\Phi_{18}=\delta_0^3+\delta_1^3+\cdots+\delta_{12}^3+\delta_{\infty}^3.$$
The calculation of $\Phi_{18}$ is similar as that of $\Phi_{12}$. We find
that
$$\Phi_{18}(x_1(z), \ldots, x_6(z))=13 \cdot 6 \Delta(z) E_6(z).$$

  Let
$$\Phi_{30}=\delta_0^5+\delta_1^5+\cdots+\delta_{12}^5+\delta_{\infty}^5.$$
As a polynomial in six variables, $\Phi_{30}(z_1, z_2, z_3, z_4, z_5, z_6)$
is a $G$-invariant polynomial. Similarly as above, we can show that
$\Phi_{30}(x_1(z), \ldots, x_6(z))$ is a modular form of weight $30$ for the
full modular group $\Gamma(1)$. Moreover, we will show that it is a cusp form.
In fact,
$$\aligned
  &\Phi_{30}(a_1(z), \ldots, a_6(z))=13^{10} q^{\frac{35}{4}} (1+O(q))^5+\\
  &+\sum_{\nu=0}^{12} [-13 q^{\frac{7}{4}} (1+O(q))+\\
  &+\zeta^{\nu} q^{\frac{43}{52}} (13+O(q))+\zeta^{2 \nu} q^{\frac{47}{52}} (-22+O(q))
   +\zeta^{3 \nu} q^{\frac{51}{52}} (-21+O(q))+\\
  &+\zeta^{4 \nu} q^{\frac{3}{52}} (-1+O(q))+\zeta^{5 \nu} q^{\frac{7}{52}} (2+O(q))
   +\zeta^{6 \nu} q^{\frac{11}{52}} (2+O(q))+\\
  &+\zeta^{7 \nu} q^{\frac{15}{52}} (-2+O(q))+\zeta^{8 \nu} q^{\frac{19}{52}} (-8+O(q))
   +\zeta^{9 \nu} q^{\frac{23}{52}} (6+O(q))+\\
  &+\zeta^{10 \nu} q^{\frac{27}{52}} (1+O(q))+\zeta^{11 \nu} q^{\frac{31}{52}} (-8+O(q))
   +\zeta^{12 \nu} q^{\frac{35}{52}} (17+O(q))]^5.
\endaligned$$
We will calculate the $q^{\frac{3}{4}}$-term which is the lowest degree.
(1) For the partition $39=4 \cdot 3+27$, the corresponding term is
$$\begin{pmatrix} 5\\ 4, 1 \end{pmatrix} (\zeta^{4 \nu} q^{\frac{3}{52}} \cdot (-1))^4 \cdot
  \zeta^{10 \nu} q^{\frac{27}{52}}=5 q^{\frac{3}{4}}.$$
(2) For the partition $39=3 \cdot 3+7+23$, the corresponding term is
$$\begin{pmatrix} 5\\ 3, 1, 1 \end{pmatrix} (\zeta^{4 \nu} q^{\frac{3}{52}} \cdot (-1))^3 \cdot
  \zeta^{5 \nu} q^{\frac{7}{52}} \cdot 2 \cdot \zeta^{9 \nu} q^{\frac{23}{52}} \cdot 6=-240 q^{\frac{3}{4}}.$$
(3) For the partition $39=3 \cdot 3+11+19$, the corresponding term is
$$\begin{pmatrix} 5\\ 3, 1, 1 \end{pmatrix} (\zeta^{4 \nu} q^{\frac{3}{52}} \cdot (-1))^3 \cdot
  \zeta^{6 \nu} q^{\frac{11}{52}} \cdot 2 \cdot \zeta^{8 \nu} q^{\frac{19}{52}} \cdot (-8)=320 q^{\frac{3}{4}}.$$
(4) For the partition $39=3 \cdot 3+2 \cdot 15$, the corresponding term is
$$\begin{pmatrix} 5\\ 3, 2 \end{pmatrix} (\zeta^{4 \nu} q^{\frac{3}{52}} \cdot (-1))^3 \cdot
  (\zeta^{7 \nu} q^{\frac{15}{52}} \cdot (-2))^2=-40 q^{\frac{3}{4}}.$$
(5) For the partition $39=2 \cdot 3+3 \cdot 11$, the corresponding term is
$$\begin{pmatrix} 5\\ 2, 3 \end{pmatrix} (\zeta^{4 \nu} q^{\frac{3}{52}} \cdot (-1))^2 \cdot
  (\zeta^{6 \nu} q^{\frac{11}{52}} \cdot 2)^3=80 q^{\frac{3}{4}}.$$
(6) For the partition $39=2 \cdot 3+2 \cdot 7+19$, the corresponding term is
$$\begin{pmatrix} 5\\ 2, 2, 1 \end{pmatrix} (\zeta^{4 \nu} q^{\frac{3}{52}} \cdot (-1))^2 \cdot
  (\zeta^{5 \nu} q^{\frac{7}{52}} \cdot 2)^2 \cdot \zeta^{8 \nu} q^{\frac{19}{52}} \cdot (-8)=-960 q^{\frac{3}{4}}.$$
(7) For the partition $39=2 \cdot 3+7+11+15$, the corresponding term is
$$\begin{pmatrix} 5\\ 2, 1, 1, 1 \end{pmatrix} (\zeta^{4 \nu} q^{\frac{3}{52}} \cdot (-1))^2 \cdot
  \zeta^{5 \nu} q^{\frac{7}{52}} \cdot 2 \cdot \zeta^{6 \nu} q^{\frac{11}{52}} \cdot 2 \cdot \zeta^{7 \nu}
  q^{\frac{15}{52}} \cdot (-2)=-480 q^{\frac{3}{4}}.$$
(8) For the partition $39=1 \cdot 3+3 \cdot 7+15$, the corresponding term is
$$\begin{pmatrix} 5\\ 1, 3, 1 \end{pmatrix} \zeta^{4 \nu} q^{\frac{3}{52}} \cdot (-1) \cdot
  (\zeta^{5 \nu} q^{\frac{7}{52}} \cdot 2)^3 \cdot \zeta^{7 \nu} q^{\frac{15}{52}} \cdot (-2)=320 q^{\frac{3}{4}}.$$
(9) For the partition $39=1 \cdot 3+2 \cdot 7+2 \cdot 11$, the corresponding term is
$$\begin{pmatrix} 5\\ 1, 2, 2 \end{pmatrix} \zeta^{4 \nu} q^{\frac{3}{52}} \cdot (-1) \cdot
  (\zeta^{5 \nu} q^{\frac{7}{52}} \cdot 2)^2 \cdot (\zeta^{6 \nu} q^{\frac{11}{52}} \cdot 2)^2=-480 q^{\frac{3}{4}}.$$
(10) For the partition $39=4 \cdot 7+11$, the corresponding term is
$$\begin{pmatrix} 5\\ 4, 1 \end{pmatrix} (\zeta^{5 \nu} q^{\frac{7}{52}} \cdot 2)^4 \cdot
  \zeta^{6 \nu} q^{\frac{11}{52}} \cdot 2=160 q^{\frac{3}{4}}.$$
Hence, for $\Phi_{30}(x_1(z), \ldots, x_6(z))$ which is a modular form for $\Gamma(1)$ with weight $30$,
the lowest degree term is given by
$$(5-240+320-40+80-960-480+320-480+160) q^{\frac{3}{4}} \cdot q^{\frac{30}{24}}=-1315 q^2.$$
Thus,
$$\Phi_{30}(x_1(z), \ldots, x_6(z))=q^2 (-13 \cdot 1315+O(q)).$$
The leading term of $\Phi_{30}(x_1(z), \ldots, x_6(z))$ together with its weight $30$ suffice to identify
this modular form with
$$\Phi_{30}(x_1(z), \ldots, x_6(z))=-13 \cdot 1315 \Delta(z)^2 E_6(z).$$

  Up to a constant, we revise the definition of $\Phi_{12}$, $\Phi_{12}^{\prime}$,
$\Phi_{18}$, $\Phi_{20}$ and $\Phi_{30}$ as given by (1.8), (1.9) and (1.10).
Consequently,
$$\left\{\aligned
  \Phi_{12}(x_1(z), \ldots, x_6(z)) &=\Delta(z),\\
  \Phi_{12}^{\prime}(x_1(z), \ldots, x_6(z)) &=\Delta(z),\\
  \Phi_{18}(x_1(z), \ldots, x_6(z)) &=\Delta(z) E_6(z),\\
  \Phi_{20}(x_1(z), \ldots, x_6(z)) &=\eta(z)^8 \Delta(z) E_4(z),\\
  \Phi_{30}(x_1(z), \ldots, x_6(z)) &=\Delta(z)^2 E_6(z).
\endaligned\right.\eqno{(3.17)}$$
This completes the proof of Theorem 1.1.

\noindent
$\qquad \qquad \qquad \qquad \qquad \qquad \qquad \qquad \qquad
 \qquad \qquad \qquad \qquad \qquad \qquad \qquad \boxed{}$

\begin{center}
{\large\bf 4. A different construction: from the modular curve
              $X(13)$ to $E_8$}
\end{center}

\noindent{\bf 4.1. $E_8$-singularity: from $X(13)$ to $E_8$}

  In this section, we will give a different construction of the
equation of the $E_8$-singularity: the symmetry group is the simple
group $\text{PSL}(2, 13)$ and the equation has a modular
interpretation in terms of theta constants of order thirteen.

  Put
$$\Phi_j=\Phi_j(x_1(z), \ldots, x_6(z)) \quad \text{for $j=12, 18, 20, 30$}.$$
By Theorem 1.1, the relations
$$j(z):=\frac{E_4(z)^3}{\Delta(z)}=\frac{\Phi_{20}^3}{\Phi_{12}^5}, \quad
  j(z)-1728=\frac{E_6(z)^2}{\Delta(z)}=\frac{\Phi_{30}^2}{\Phi_{12}^5}\eqno{(4.1)}$$
give the equation
$$\Phi_{20}^3-\Phi_{30}^2=1728 \Phi_{12}^5.\eqno{(4.2)}$$
Hence, we have the following:

\vskip 0.1 cm

\textbf{Theorem 4.1} (A different construction of the $E_8$-singularity:
from $X(13)$ to $E_8$). {\it The equation of the $E_8$-singularity can
be constructed from the modular curve $X(13)$ as follows$:$
$$\Phi_{20}^3-\Phi_{30}^2=1728 \Phi_{12}^5,$$
where $\Phi_{12}$, $\Phi_{20}$ and $\Phi_{30}$ are $G$-invariant
polynomials.}

  Let us recall some facts about exotic spheres (see \cite{Hi}).
A $k$-dimensional compact oriented differentiable manifold is
called a $k$-sphere if it is homeomorphic to the $k$-dimensional
standard sphere. A $k$-sphere not diffeomorphic to the standard
$k$-sphere is said to be exotic. The first exotic sphere was
discovered by Milnor in 1956 (see \cite{Mi}). Two
$k$-spheres are called equivalent if there exists an orientation
preserving diffeomorphism between them. The equivalence classes
of $k$-spheres constitute for $k \geq 5$ a finite abelian group
$\Theta_k$ under the connected sum operation. $\Theta_k$ contains
the subgroup $b P_{k+1}$ of those $k$-spheres which bound a
parallelizable manifold. $b P_{4m}$ ($m \geq 2$) is cyclic of order
$2^{2m-2}(2^{2m-1}-1)$ numerator $(4 B_m/m)$, where $B_m$ is the
$m$-th Bernoulli number. Let $g_m$ be the Milnor generator of $b P_{4m}$.
If a $(4m-1)$-sphere $\Sigma$ bounds a parallelizable manifold $B$ of
dimension $4m$, then the signature $\tau(B)$ of the intersection form
of $B$ is divisible by $8$ and $\Sigma=\frac{\tau(B)}{8} g_m$. For $m=2$
we have $b P_8=\Theta_7=\mathbb{Z} /28 \mathbb{Z}$. All these results are
due to Milnor-Kervaire (see \cite{KM}). In particular,
$$\sum_{i=0}^{2m} z_i \overline{z_i}=1, \quad
  z_0^3+z_1^{6k-1}+z_2^2+\cdots+z_{2m}^2=0$$
is a $(4m-1)$-sphere embedded in $S^{4m+1} \subset \mathbb{C}^{2n+1}$
which represents the element $(-1)^m k \cdot g_m \in b P_{4m}$. For
$m=2$ and $k=1, 2, \cdots, 28$ we get the $28$ classes of $7$-spheres.
Theorem 4.1 shows that the higher dimensional liftings of two distinct
symmetry groups and modular interpretations on the equation of the
$E_8$-singularity give the same Milnor's standard generator of $\Theta_7$.

\vskip 0.3 cm

\noindent{\bf 4.2. The Barlow surface: from $X(13)$ to $E_8$}

  In this section, we will construct the Barlow surface from the
modular curve $X(13)$ by the method of transversal linear sections.

  By Theorem 1.1, we have
$$\Phi_{12}(x_1(z), \ldots, x_6(z)) \cdot \Phi_{18}(x_1(z), \ldots, x_6(z))
 -\Phi_{30}(x_1(z), \ldots, x_6(z))=0.\eqno{(4.3)}$$
This gives a morphism
$$\Phi: X(13) \to X \subset \mathbb{CP}^5\eqno{(4.4)}$$
with $\Phi(z)=(x_1(z), \ldots, x_6(z))$, where the variety $X$ is
given by an equation in six variables $z_1, \ldots, z_6$ of degree
$30$:
$$\Phi_{12}(z_1, \ldots, z_6) \cdot \Phi_{18}(z_1, \ldots, z_6)
 -\Phi_{30}(z_1, \ldots, z_6)=0.\eqno{(4.5)}$$
Moreover, we have the following morphism
$$\aligned
  \varphi: X &\to Z \subset \mathbb{CP}^{13},\\
           (z_1, \ldots, z_6) &\mapsto
           (\delta_0, \delta_1, \ldots, \delta_{12}, \delta_{\infty}),
\endaligned\eqno{(4.6)}$$
where $\delta_0, \delta_1, \ldots, \delta_{12}, \delta_{\infty}$ are
given by (1.4). The variety $Z$ is given by
$$\left\{\aligned
  \delta_0+\delta_1+\cdots+\delta_{12}+\delta_{\infty} &=0,\\
  (\delta_0^2+\cdots+\delta_{\infty}^2)(\delta_0^3+\cdots+\delta_{\infty}^3)
 -\frac{4056}{1315}(\delta_0^5+\cdots+\delta_{\infty}^5) &=0,
\endaligned\right.\eqno{(4.7)}$$
which is a Fano $11$-fold. Let
$$W=Z \cap H_1 \cap H_2 \cap \cdots \cap H_9\eqno{(4.8)}$$
be a transversal linear section, where the hyperplane
$H_i:=\{ \delta_{4+i}=0 \}$ $(1 \leq i \leq 8)$ and
$H_9:=\{ \delta_{\infty}=0 \}$. Now, let us recall the following:

\textbf{Proposition 4.2} (see \cite{GZ}). {\it In the family of quintic
surfaces $S_{(\lambda: \mu)}$ in $\mathbb{P}^4$ given by
$$S_{(\lambda: \mu)}=\{ (x_0: \ldots :x_4) \in \mathbb{P}^4:
  \sigma_1=0, \lambda \sigma_2 \sigma_3+\mu \sigma_5=0 \},$$
$($where $\sigma_k$ denotes the $k$-th elementary symmetric
polynomial in $x_0, \ldots, x_4$ and $\lambda, \mu \in \mathbb{C}$$)$,
all but the following six are non-singular$:$

$\mathrm{(i)}$ $\sigma_5=0$ $($reducible, consisting of $5$ planes, meeting
along $10$ lines which in turn meet $3$ at a time in $10$ points$)$,

$\mathrm{(ii)}$ $\sigma_2 \sigma_3=0$ $($reducible, consisting of a quadric
and a cubic surface meeting along a non-singular sextic curve$)$,

$\mathrm{(iii)}$ $2 \sigma_5+\sigma_2 \sigma_3=0$ $($$20$ singularities, namely
the $S_5$-orbit of $(-2: -2: -2: 3+\sqrt{-7}: 3-\sqrt{-7})$$)$,

$\mathrm{(iv)}$ $25 \sigma_5-12 \sigma_2 \sigma_3=0$ $($$10$ singularities,
namely the $S_5$-orbit of $(-2: -2: -2: 3 : 3)$$)$,

$\mathrm{(v)}$ $50 \sigma_5+\sigma_2 \sigma_3=0$ $($$5$ singularities, namely
the $S_5$-orbit of $(1 : 1 : 1 : 1 : -4)$$)$,

$\mathrm{(vi)}$ $2 \sigma_5-\sigma_2 \sigma_3=0$ $($$15$ singularities, namely
the $S_5$-orbit of $(0 : 1 : -1 : 1 : -1)$$)$.}

  Using the notation in Theorem 4.2, we find that $W$ is equivalent to
the quintic surface $S_{(1: -\frac{676}{413})}$, which can be deformed
in the family $S_{(\lambda: \mu)}$ to the surface (iii) in Theorem 4.2.
Note that the surface (iii) in Theorem 4.2 is equivalent to the
quintic surface $Q$ in (2.5). Let $\Phi: Y \to Q$ be a covering map,
where $Y$ is $2:1$ onto the $20$-nodal quintic $Q \subset \mathbb{CP}^4$.
The icosahedral group $A_5$ acts on $\mathbb{CP}^4$ by the standard
action on the coordinates. The quintic $Q$ is $A_5$-invariant and
its $20$ nodes are the $A_5$-orbit of the point $(2, 2, 2,
-3-\sqrt{-7}, -3+\sqrt{-7})$. The $A_5$-action on $Q$ is covered
by an action on $Y$, so that we have an action of $A_5 \times
\mathbb{Z}/2 \mathbb{Z}=A_5 \cup A_5 \sigma$ on $Y$, where the
generator $\sigma \in \mathbb{Z}/2 \mathbb{Z}$ is the covering
involution. Elements of $A_5 \times \mathbb{Z}/2 \mathbb{Z}$
acting on $Y$ are denoted like the corresponding elements acting
on $Q$. Now, we can apply the main result from \cite{Ba}:

\textbf{Proposition 4.3} (see \cite{Ba} and \cite{Ko}). {\it Let
$\alpha=(02)(34) \sigma$, $\beta=(01234)$. Then $\beta$ acts freely
on $Y$ and $\alpha$ has $4$ fixed points. The resolution of the
nodes of $Y/D_{10}$, where $D_{10}=\langle \alpha, \beta \rangle$,
gives a minimal surface $B$ of general type with $\pi_1=0$, $q=p_g=0$
and $K^2=1$.}

  By \cite{Ko}, $B$ is homeomorphic but not diffeomorphic to
$\mathbb{CP}^2 \# 8 \overline{\mathbb{CP}}^2$. Therefore, we have
proved the following:

\textbf{Theorem 4.4.} {\it The Barlow surface can be constructed
from the modular curve $X(13)$.}

\vskip 0.3 cm

\noindent{\bf 4.3. Bring's curve and Fricke's octavic curve: from
                $X(13)$ to $E_8$}

  In this section, we will construct both Bring's curve and Fricke's
octavic curve from the modular curve $X(13)$ by the method of
transversal linear sections.

  By Theorem 1.1, we have
$$\left\{\aligned
  \Phi_{4}(x_1(z), \ldots, x_6(z)) &=0,\\
  \Phi_{8}(x_1(z), \ldots, x_6(z)) &=0,\\
  \phi_{12}(x_1(z), \ldots, x_6(z)) &=0,
\endaligned\right.$$
which gives a morphism from $X(13)$ to an algebraic surface $S_1
\subset \mathbb{CP}^5$ given by the following $G$-invariant equations
$$\left\{\aligned
  \Phi_{4}(z_1, \ldots, z_6) &=0,\\
  \Phi_{8}(z_1, \ldots, z_6) &=0,\\
  \phi_{12}(z_1, \ldots, z_6) &=0.
\endaligned\right.$$
Moreover, there is a morphism
$$\aligned
  \varphi_1: S_1 &\to Y_1 \subset \mathbb{CP}^{27},\\
           (z_1, \ldots, z_6) &\mapsto (w_0, w_1, \ldots, w_{12},
           w_{\infty}, \delta_0, \delta_1, \ldots, \delta_{12},
           \delta_{\infty}),
\endaligned\eqno{(4.9)}$$
where $w_0, w_1, \ldots, w_{12}, w_{\infty}$ are given by (1.2) and
$\delta_0, \delta_1, \ldots, \delta_{12}, \delta_{\infty}$ are given
by (1.4), and the variety $Y_1$ is given by
$$\left\{\aligned
  w_0+w_1+\cdots+w_{12}+w_{\infty} &=0,\\
  w_0^2+w_1^2+\cdots+w_{12}^2+w_{\infty}^2 &=0,\\
  w_0^3+w_1^3+\cdots+w_{12}^3+w_{\infty}^3 &=\frac{15}{26}
  (\delta_0^2+\delta_1^2+\cdots+\delta_{12}^2+\delta_{\infty}^2),
\endaligned\right.\eqno{(4.10)}$$
which is a $24$-dimensional Fano variety. Let
$$C_1=Y_1 \cap H_1 \cap H_2 \cap \cdots \cap H_{23}\eqno{(4.11)}$$
be a transversal linear section, where the hyperplane
$H_i:=\{ w_{4+i}=0 \}$ $(1 \leq i \leq 8)$,
$H_9:=\{ w_{\infty}=0 \}$, $H_j:=\{ \delta_{j-10}=0 \}$
$(10 \leq j \leq 22)$ and $H_{23}:=\{ \delta_{\infty}=0 \}$.
Then $C_1$ is equivalent to Bring's curve.

  Similarly, by Theorem 1.1, we have
$$\left\{\aligned
  \Phi_{4}(x_1(z), \ldots, x_6(z)) &=0,\\
  \Phi_{8}(x_1(z), \ldots, x_6(z)) &=0,\\
  \Phi_{16}(x_1(z), \ldots, x_6(z)) &=0,
\endaligned\right.$$
which gives a morphism from $X(13)$ to an algebraic surface $S_2
\subset \mathbb{CP}^5$ given by the following $G$-invariant equations
$$\left\{\aligned
  \Phi_{4}(z_1, \ldots, z_6) &=0,\\
  \Phi_{8}(z_1, \ldots, z_6) &=0,\\
  \Phi_{16}(z_1, \ldots, z_6) &=0.
\endaligned\right.$$
Moreover, there is a morphism
$$\aligned
  \varphi_2: S_2 &\to Y_2 \subset \mathbb{CP}^{13},\\
           (z_1, \ldots, z_6) &\mapsto (w_0, w_1, \ldots, w_{12}, w_{\infty}),
\endaligned\eqno{(4.12)}$$
where $w_0, w_1, \ldots, w_{12}, w_{\infty}$ are given by (1.2) and the
variety $Y_2$ is given by
$$\left\{\aligned
  w_0+w_1+\cdots+w_{12}+w_{\infty} &=0,\\
  w_0^2+w_1^2+\cdots+w_{12}^2+w_{\infty}^2 &=0,\\
  w_0^4+w_1^4+\cdots+w_{12}^4+w_{\infty}^4 &=0,
\endaligned\right.\eqno{(4.13)}$$
which is a Fano $10$-fold. Let
$$C_2=Y_2 \cap H_1 \cap H_2 \cap \cdots \cap H_9\eqno{(4.14)}$$
be a transversal linear section, where the hyperplane
$H_i:=\{ w_{4+i}=0 \}$ $(1 \leq i \leq 8)$ and
$H_9:=\{ w_{\infty}=0 \}$. Then $C_2$ is equivalent
to Fricke's octavic curve. Therefore, we have proved
the following:

\textbf{Theorem 4.5.} {\it Both Bring's curve and Fricke's
octavic curve can be constructed from the modular curve $X(13)$.}

\begin{center}
{\large\bf 5. An explicit construction of the modular curve $X(13)$}
\end{center}

  In this section, we will study the following classical problem in
the theory of modular curves, which goes back to F. Klein (see \cite{K1}
and \cite{K2}).

\textbf{Problem 5.1.} Let $p \geq 7$ be a prime number. Give an
explicit construction of the modular curve $X(p)$ of level $p$
from the invariant theory for $\mathrm{PSL}(2, p)$ using projective
algebraic geometry.

\textbf{Example 5.2.} The equation of the modular curve $X(7)$ of
level $7$ is given by the Klein quartic curve (see \cite{K1})
$$x^3 y+y^3 z+z^3 x=0\eqno{(5.1)}$$
in $\mathbb{P}^2$. The left hand side of (5.1) is the unique quartic
invariant for $\mathrm{PSL}(2, 7)$ in this representation.

\textbf{Example 5.3.} The matrix
$$\begin{pmatrix}
  w & v & 0 & 0 & z\\
  v & x & w & 0 & 0\\
  0 & w & y & x & 0\\
  0 & 0 & x & z & y\\
  z & 0 & 0 & y & v
  \end{pmatrix}$$
is (up to a factor) the Hessian matrix of a cubic invariant for
$\mathrm{PSL}(2, 11)$, namely the Klein cubic threefold (see \cite{K2})
$$v^2 w+w^2 x+x^2 y+y^2 z+z^2 v=0.\eqno{(5.2)}$$
The modular curve $X(11)$ of level $11$ is the singular locus of the
Hessian of this cubic threefold. In fact, there are only $10$ distinct
quartics defining the locus, namely
$$v^2 wz-v w y^2-x^2 yz=0, \quad vy^3+w^3 z+w x^3=0$$
and their images under successive applications of the cyclic permutation
$(vwxyz)$.

  In general, the locus of the modular curve $X(p)$ of level $p$ can
be defined by $\begin{pmatrix} (p-1)/2 \\ 3 \end{pmatrix}$ quartics
(see \cite{AR}, p. 59). Following Klein's method for the cubic threefold
(5.2), Adler and Ramanan (see \cite{AR}) studied Problem 5.1 when $p$ is a
prime congruent to $3$ modulo $8$ by some cubic hypersurface invariant
under $\mathrm{PSL}(2, p)$. However, their method can not be valid
for $p=13$. In this case, the number of the locus is $20$. As Ramanan
has remarked, one can actually reduce the number further, namely to
$(p-1)/2$, but this does not lead to explicit equations (see \cite{AR},
p. 59).

  As a consequence of Theorem 1.1, we find an explicit construction
of the modular curve $X(13)$.

\textbf{Theorem 5.4.} {\it There is a morphism
$$\Phi: X(13) \to C \subset \mathbb{CP}^5\eqno{(5.3)}$$
with $\Phi(z)=(x_1(z), \ldots, x_6(z))$, where $C$ is an algebraic
curve given by a family of $G$-invariant equations
$$\left\{\aligned
  \Phi_{4}(z_1, \ldots, z_6) &=0,\\
  \Phi_{8}(z_1, \ldots, z_6) &=0,\\
  \phi_{12}(z_1, \ldots, z_6) &=0,\\
  \Phi_{16}(z_1, \ldots, z_6) &=0.
\endaligned\right.\eqno{(5.4)}$$}

{\it Proof}. Theorem 1.1 implies that
$$\left\{\aligned
  \Phi_{4}(x_1(z), \ldots, x_6(z)) &=0,\\
  \Phi_{8}(x_1(z), \ldots, x_6(z)) &=0,\\
  \phi_{12}(x_1(z), \ldots, x_6(z)) &=0,\\
  \Phi_{16}(x_1(z), \ldots, x_6(z)) &=0.
\endaligned\right.$$
\noindent
$\qquad \qquad \qquad \qquad \qquad \qquad \qquad \qquad \qquad
 \qquad \qquad \qquad \qquad \qquad \qquad \qquad \boxed{}$

\vskip 2.0 cm

\noindent{Department of Mathematics, Peking University}

\noindent{Beijing 100871, P. R. China}

\noindent{\it E-mail address}: yanglei$@$math.pku.edu.cn

\end{document}